\documentclass[preprint,12pt]{amsart}

\usepackage[utf8]{inputenc}  
\usepackage[T1]{fontenc}     
\usepackage{tikz}
\usetikzlibrary{arrows.meta, positioning}

\usepackage[foot]{amsaddr}
\usepackage[margin=1in]{geometry}

\usepackage{amsmath, amsfonts, amsthm, bbm}
\usepackage{enumerate}
\numberwithin{equation}{section}

\usepackage{graphicx}
\graphicspath{{./images/}{./figures/}}
\usepackage{subcaption}
\usepackage{tikz}
\usepackage{url}
\usepackage{algorithmic}

\usepackage[linesnumbered,ruled,vlined]{algorithm2e}
\SetKwInput{KwInput}{Input}
\SetKwInput{KwOutput}{Output}

\usepackage{xcolor}

\newtheorem{thm}{Theorem}[section]

\newtheorem{rem}[thm]{Remark}

\def\XXint#1#2#3{{\setbox0=\hbox{$#1{#2#3}{\int}$}
     \vcenter{\hbox{$#2#3$}}\kern-.5\wd0}}

\title[Hybrid PINNs for Chemotaxis]{A Hybrid Physics-Informed Neural Network Framework for Subcritical and Supercritical Dynamics in Multi-Species Chemotaxis}

\author{Hamid El Bahja}
\address{African Institute for Mathematical Sciences, Research and Innovation Centre, Kigali, Rwanda}
\email{helbahja@aims.ac.za}

\author{Jan C. Riedel}
\address{Technical University Berlin, Germany}
\email{j.hauffen@tu-berlin.de}

\author{Peter Jung}
\address{Technical University Berlin \& German Aerospace Center (DLR), Germany}
\email{peter.jung@tu-berlin.de}

\begin{document}

\begin{abstract}
We study a two-species chemotaxis system that exhibits two qualitatively different regimes: subcritical dynamics, in which solutions remain smooth, and supercritical dynamics, in which strong aggregation may lead to finite-time blow-up. In the subcritical regime, we use a standard continuous-time physics-informed neural network (PINN) with alternating training for the coupled fields and show that it provides accurate and efficient approximations. In the supercritical regime, however, this formulation is not sufficiently robust to resolve the highly localized structures and steep gradients associated with blow-up. To address this limitation, we introduce a discrete-time PINN based on backward Euler time stepping, combined with a logarithmic transformation to stabilize large solution values and a residual-adaptive collocation strategy that concentrates training points near regions of high activity. This hybrid framework allows us to use a simple and computationally inexpensive PINN when the solution is smooth, while switching to a more robust discrete formulation when singular dynamics emerge. Numerical experiments confirm that the continuous-time formulation accurately resolves subcritical dynamics, while the discrete-time formulation reliably captures singular aggregation and finite-time blow-up in the supercritical regime.
\end{abstract}

\keywords{Physics-informed neural networks, chemotaxis, blow-up phenomena, backward Euler method}

\maketitle

\section{Introduction}

Chemotaxis is the process by which motile cells sense and migrate along chemical gradients in their surrounding environment. By detecting variations in attractant or repellent concentrations, cells bias their movement toward favorable regions (e.g., nutrient-rich areas or signaling sources) while avoiding harm. This mechanism underpins many biological processes, including bacterial aggregation, immune responses, embryonic development, wound healing, and tumor invasion~\cite{Devreotes1989,Murray2002}. In many real-world settings, however, multiple species, such as competing bacterial strains, different immune cell types, or heterogeneous tumor subpopulations, respond to the same chemical signal. This naturally leads to multi-species chemotaxis models, where distinct cell densities are coupled through a common chemical field, giving rise to richer dynamics like cooperative aggregation, competition, and spatial segregation~\cite{Horstmann2011,EspejoArenas2009,Ucar2025,Kolbe2020}.

The mathematical study of chemotaxis dates back to Patlak~\cite{Patlak1953}, who described biased random walks at the macroscopic level. Building on this idea, Keller and Segel introduced a system of partial differential equations in the 1970s to model the coupled evolution of cell density and chemoattractant concentration, capturing the feedback loop between cell movement and signal dynamics~\cite{KellerSegel1971}. Since then, Keller–Segel type models have become a cornerstone of mathematical biology, extensively studied for their rich analytical behavior, including pattern formation, aggregation, and finite-time blow-up~\cite{Childress1981,HerreroVelazquez1997,Hillen2009,Perthame2006, ElBahja2024}.

In this work, we consider a two-species chemotaxis model in which two interacting populations are coupled to a single chemoattractant field. The evolution equations take the form
\begin{equation}\label{eq:full_system_intro}
\begin{cases}
\partial_t \rho_1 = \mu_1 \Delta \rho_1 - \chi_1 \nabla\cdot(\rho_1 \nabla c),\\[4pt]
\partial_t \rho_2 = \mu_2 \Delta \rho_2 - \chi_2 \nabla\cdot(\rho_2 \nabla c),\\[4pt]
\varepsilon\,\partial_t c = D\Delta c + \alpha_1 \rho_1 + \alpha_2 \rho_2 - \beta c,
\end{cases}
\end{equation}
for \(x\in\Omega\subset\mathbb{R}^d\) and \(t>0\), where \(\rho_1\) and \(\rho_2\) denote the densities of the two populations and \(c\) the chemoattractant concentration. The parameters \(\mu_i>0\) represent diffusion coefficients, \(\chi_i>0\) the chemotactic sensitivities, while \(\alpha_i\), \(D\), and \(\beta\) denote the production, diffusion, and degradation rates of the chemical signal.

The system is complemented with homogeneous Neumann boundary conditions
\begin{equation}
\nabla \rho_i \cdot \nu = 0, \qquad \nabla c \cdot \nu = 0
\quad \text{on } \partial\Omega ,
\end{equation}
together with initial data
\begin{equation}
\rho_i(0,x)=\rho_{i,0}(x)~~\text{for}~i=1,2, \qquad c(0,x)=c_0(x).
\end{equation}
The no–flux conditions ensure that neither the cell populations nor the chemoattractant leave the domain.

The parameter \(\varepsilon\geq0\) controls the time scale of the chemoattractant dynamics. For \(\varepsilon>0\), system \eqref{eq:full_system_intro} is fully parabolic. In the limit \(\varepsilon=0\), the chemoattractant equilibrates instantaneously and the third equation reduces to the elliptic relation
\[
-D\Delta c = \alpha_1\rho_1+\alpha_2\rho_2-\beta c,
\]
which yields a parabolic–elliptic chemotaxis system. Moreover, if one species is absent, for example \(\rho_2\equiv0\), the model reduces to the classical parabolic chemotaxis system introduced by Patlak and Keller–Segel~\cite{Patlak1953,KellerSegel1971}. Hence, \eqref{eq:full_system_intro} provides a unified framework encompassing both multi-species and classical chemotaxis dynamics.

Because explicit analytical solutions are unavailable for most multi‑species and blow‑up scenarios, numerical simulation is essential. Depending on the parameter regime and initial data, chemotaxis systems may develop strongly localized aggregation peaks and, in some cases, finite‑time blow‑up~\cite{Winkler2010,Horstmann2011}. Capturing these singular dynamics is numerically challenging: the solution develops steep gradients and very large amplitudes, while the parabolic‑elliptic limit and multi‑species coupling introduce additional stiffness.

Traditional numerical approaches for chemotaxis models include finite element, finite volume, and operator‑splitting methods for single‑species systems~\cite{xiao2019numerical,li2017local,guo2019energy,chertock2008second,filbet2006finite,tyson2000fractional}. For multi‑species models, the two‑species parabolic‑elliptic and parabolic‑parabolic systems have been studied numerically using high‑order finite‑volume and adaptive moving mesh schemes~\cite{KurganovLukacovaMedvidova2014,Chertock2019}. These methods can capture different blow‑up scenarios, but they require very fine uniform meshes (up to $400\times400$) to resolve the disparate growth rates of the two species, leading to high computational cost. 

To overcome these limitations, deep learning‑based approaches have recently emerged as promising alternatives for solving partial differential equations. In particular, physics‑informed neural networks (PINNs) incorporate the governing equations directly into the loss function, enabling mesh‑free approximation through automatic differentiation~\cite{rag}. These methods have been successfully applied to a wide range of nonlinear and high‑dimensional PDEs~\cite{han, ElBahja2025}.

However, standard PINNs often struggle in the presence of sharp gradients, steep fronts, and localized singularities, where the underlying solution exhibits rapid spatial or temporal variations that are difficult to capture accurately with smooth neural network approximations. To address these limitations, several improvements have been proposed, including residual‑based adaptive sampling strategies that concentrate collocation points in regions of large error~\cite{Wu2023,MaoMeng2023,Gao2023}, as well as output transformations and normalization techniques designed to stabilise training in the presence of large solution values~\cite{LeDuc2024}. Despite these advances, accurately resolving strongly localised aggregation and multi‑scale dynamics in chemotaxis systems remains challenging.

Motivated by this observation, we propose a hybrid PINN framework that adapts to the qualitative behavior of the solution. The main idea is to use the simplest effective neural formulation in each regime rather than relying on one uniform method for all cases. More precisely, our approach provides the following advantages:
\begin{itemize}
    \item Efficiency in smooth regimes: for subcritical dynamics, a standard continuous-time PINN with alternating training is employed, offering computational efficiency while accurately capturing regular solutions.
    
    \item Robustness near blow-up: for supercritical dynamics, we introduce a discrete-time PINN based on backward Euler time stepping, which provides superior stability for stiff problems with rapidly evolving solutions.
    
    \item Improved numerical stability: a logarithmic transformation enforces solution positivity, compresses the dynamic range of large amplitudes, and facilitates the representation of strongly localized peaks by the neural network.
    
    \item Targeted resolution of singular structures: residual-adaptive collocation focuses training points in regions of high activity, ensuring that the network allocates its capacity where it is most needed to resolve steep gradients and incipient singularities.
    
    \item Unified treatment of multi-species dynamics: the same framework seamlessly handles two interacting cell populations coupled through a common chemoattractant field, preserving the coupling between species while accurately resolving their disparate aggregation scales.
\end{itemize}
This hybrid design is the central contribution of the paper. Rather than forcing one PINN architecture to perform equally well in all regimes, we tailor the training strategy to the underlying dynamics. As a result, the proposed method combines the computational economy of standard PINNs in smooth settings, owing to their simpler architecture, single global training pass, and absence of time-stepping overhead, with the robustness of a discrete adaptive formulation in singular settings. Numerical experiments demonstrate that this strategy accurately captures both globally smooth solutions and strongly localized blow-up patterns in multi-species chemotaxis.

The remainder of the paper is organized as follows. In Section~2, we review known analytical results on global existence and finite-time blow-up for Keller--Segel type systems, with particular emphasis on the two-species setting. Section~3 introduces the proposed hybrid PINN framework: we first present the continuous-time formulation with alternating training for subcritical regimes, and then develop the discrete-time PINN with logarithmic transformation and adaptive collocation for blow-up dynamics. In Section~4, we validate the approach through a series of numerical experiments covering subcritical, and supercritical regimes, demonstrating the accuracy and robustness of the method. Finally, Section~5 summarizes the main findings and discusses possible directions for future work.

\section{Known results on blow-up versus global existence}
\label{blow}
We recall standard results on blow-up and global existence for Keller--Segel type systems, and then state the results that are directly relevant for this work.

\subsection{Classical Keller--Segel system in two space dimensions}

We consider the two-dimensional Keller--Segel system with a parameter $\gamma \in \{0,1\}$:
\begin{equation}
\label{eq:KS_gamma}
\left\{
\begin{aligned}
\partial_t u &= \Delta u - \nabla \cdot (u \nabla v),\\
\gamma\,\partial_t v &= \Delta v - v + u,
\end{aligned}
\right.
\end{equation}
posed on a bounded smooth domain $\Omega \subset \mathbb{R}^2$ (or on $\mathbb{R}^2$), supplemented with nonnegative initial data and homogeneous Neumann boundary conditions. The total mass
\[
M = \int_\Omega u_0
\]
is conserved in time.

In two space dimensions, the system exhibits a critical mass phenomenon, with different thresholds depending on the value of $\gamma$.

\medskip
\noindent
Parabolic--elliptic case ($\gamma = 0$): the critical mass is sharp
\begin{itemize}
\item If $M < 8\pi$, all solutions exist globally in time and remain bounded.
\item If $M > 8\pi$, finite-time blow-up occurs for a large class of initial data, for example radially symmetric and sufficiently concentrated data.
\end{itemize}
This threshold is classical and has been established in a series of works; see, for example, \cite{JagerLuckhaus1992,Nagai2001,BlanchetDolbeaultPerthame2006}.

\medskip
\noindent
Fully parabolic case ($\gamma = 1$): The situation is more delicate and only partially understood:
\begin{itemize}
\item If $M < 4\pi$, every solution exists globally in time and remains bounded \cite{NagaiSenbaYoshida1997}.
\item For any $\varepsilon > 0$, there exist solutions with mass $M < 4\pi + \varepsilon$ that become unbounded \cite{HorstmannWang2001}. It is not known whether such solutions blow up in finite time or remain global while becoming unbounded as $t \to \infty$.
\item For radially symmetric solutions with $M > 8\pi$, finite-time blow-up occurs; moreover, the asymptotic blow-up profile has been characterized \cite{HerreroVelazquez1997}.
\end{itemize}

Thus, while $8\pi$ acts as a threshold for finite-time blow-up in the radially symmetric setting, no sharp mass threshold for finite-time blow-up is known in the non-radial case. In particular, the existence of unbounded solutions above $4\pi$ does not imply finite-time blow-up.

\medskip
In higher dimensions ($n \ge 3$), the situation is different: for the fully parabolic system, finite-time blow-up can occur for arbitrarily small mass provided the initial energy is sufficiently negative \cite{Winkler2013}.

\subsection{Two-species system}

We now turn to the two-species system \eqref{eq:full_system_intro}. Let
\[
M_i := \int_{\mathbb{R}^2} \rho_{i,0}(x)\,dx, \qquad i=1,2,
\]
denote the conserved masses. We assume, without loss of generality, that
\(\chi_1 \le \chi_2\). The results below concern the whole space
\(\mathbb{R}^2\).

\medskip
\noindent
Parabolic--elliptic case (\(\varepsilon=0\)): A sharp threshold curve in the
\((M_1,M_2)\)-plane has been established:
\begin{itemize}
  \item If
  \[
  \frac{4\pi\mu_1 M_1}{\chi_1} + \frac{4\pi M_2}{\chi_2}
  > \frac12(M_1+M_2)^2,
  \]
  together with
  \[
  M_1 < \frac{8\pi\mu_1}{\chi_1}, \qquad M_2 < \frac{8\pi}{\chi_2},
  \]
  then the solution exists globally in time and remains bounded
  \cite{EspejoVilchesConca2013}.
  
  \item If
  \[
  \frac{4\pi\mu_1 M_1}{\chi_1} + \frac{4\pi M_2}{\chi_2}
  < \frac12(M_1+M_2)^2,
  \]
  then finite-time blow-up occurs, and both densities become unbounded simultaneously
  \cite{EspejoVilchesConca2013}.
  
  \item If
  \[
  \frac{4\pi\mu_1 M_1}{\chi_1} + \frac{4\pi M_2}{\chi_2}
  > \frac12(M_1+M_2)^2
  \quad\text{but}\quad
  M_2 \ge \frac{8\pi}{\chi_2}\;\; \text{(or } M_1 \ge \tfrac{8\pi\mu_1}{\chi_1}\text{)},
  \]
  the theory does not prove either global existence or finite-time blow-up.
  Numerical evidence suggests that the species exceeding its individual threshold
  may blow up while the other remains bounded or grows only algebraically
  \cite{KurganovLukacovaMedvidova2014}.
\end{itemize}
The blow-up region was first identified in \cite{ConcaEspejoVilches2011}. In the radially symmetric setting, blow-up can occur even when the total second moment is increasing, in particular when \(M_1 > \frac{8\pi\mu_1}{\chi_1}\) or \(M_2 > \frac{8\pi}{\chi_2}\) \cite{ConcaEspejoVilches2011}. Moreover, when blow-up occurs, both species blow up at the same time, and the blow-up measure satisfies a mass quantization property, with possible subcollapse phenomena \cite{EspejoStevensSuzuki2012}.

\medskip
\noindent
Parabolic--parabolic case (\(\varepsilon>0\)): The qualitative picture remains similar: a critical curve in the \((M_1,M_2)\)-plane separates global existence from finite-time blow-up. Numerical simulations suggest that the species with the larger chemotactic sensitivity typically blows up first \cite{KurganovLukacovaMedvidova2014}.

\section{Methodology}
\label{sec:methodology}

This section outlines the two neural network frameworks developed in this work. For subcritical regimes, we employ a standard continuous-time PINNs formulation. To handle supercritical dynamics and capture finite-time blow-up, we introduce a discrete-time PINNs approach enhanced with residual-adaptive sampling and a logarithmic output transformation.

\subsection{Physics-Informed Neural Network with Alternating Training}
\label{sec:continuous_pinn}
We solve the multi-species chemotaxis system \eqref{eq:full_system_intro} using a continuous-time Physics-Informed Neural Network (PINN). Since the governing equations are strongly coupled and nonlinear, we employ an alternating training strategy in which one network is updated at a time while the remaining networks are kept fixed.

Following Raissi et al.~\cite{rag}, we approximate the unknown fields $\rho_1$, $\rho_2$, and $c$ by three feed-forward neural networks of depth $L$, each consisting of one input layer, $L-1$ hidden layers, and one output layer. For a generic scalar field $u(x,t)$, the neural network approximation is written as
\begin{equation}
u(x,t;\theta)
=
\left(\mathcal{L}_L \circ \sigma^{L-1} \circ \mathcal{L}_{L-1} \circ \cdots \circ \sigma^1 \circ \mathcal{L}_1\right)(x,t),
\end{equation}
where $\circ$ denotes function composition, so that $(f \circ g)(x)=f(g(x))$. Each layer is an affine map of the form
\begin{equation}
\mathcal{L}_i(z)=W^i z+b^i, \qquad i=1,\dots,L,
\end{equation}
and $\sigma^i$ denotes the activation function applied in the hidden layers. The trainable parameters are collected in
\[
\theta=\{W^i,b^i\}_{i=1}^L.
\]
In all experiments, we use the hyperbolic tangent activation function in the hidden layers, since its smoothness is well suited for automatic differentiation. The output layer is linear, and a non-negative activation map, such as the softplus function, is applied to the outputs to ensure positivity of $\rho_1$, $\rho_2$, and $c$, which is required by the model.

We denote the neural approximations by
\begin{equation}
\mathcal{N}_{\rho_1}(x,t;\theta_{\rho_1}), \qquad
\mathcal{N}_{\rho_2}(x,t;\theta_{\rho_2}), \qquad
\mathcal{N}_{c}(x,t;\theta_c).
\end{equation}
Substituting these approximations into the PDE system and differentiating with automatic differentiation yields the residuals
\begin{align}
r_{\rho_1}(x,t;\theta_{\rho_1},\theta_c)
&=
\partial_t \mathcal{N}_{\rho_1}
-\mu_1 \Delta \mathcal{N}_{\rho_1}
+\chi_1 \nabla \cdot \big(\mathcal{N}_{\rho_1}\nabla \mathcal{N}_c\big), \\
r_{\rho_2}(x,t;\theta_{\rho_2},\theta_c)
&=
\partial_t \mathcal{N}_{\rho_2}
-\mu_2 \Delta \mathcal{N}_{\rho_2}
+\chi_2 \nabla \cdot \big(\mathcal{N}_{\rho_2}\nabla \mathcal{N}_c\big), \\
r_c(x,t;\theta_c,\theta_{\rho_1},\theta_{\rho_2})
&=
\varepsilon\,\partial_t \mathcal{N}_c
-D\Delta \mathcal{N}_c
-\alpha_1 \mathcal{N}_{\rho_1}
-\alpha_2 \mathcal{N}_{\rho_2}
+\beta \mathcal{N}_c.
\end{align}
Each residual depends on the parameters of the field being updated as well as on the frozen parameters of the coupled networks. This structure motivates the alternating optimization strategy.

Training all fields jointly or updating three separate networks simultaneously can be unstable in strongly coupled systems. A similar difficulty was reported by Niaki et al. \cite{AminiNiaki2021}, where joint training failed to sufficiently reduce the loss of one variable and, in the case of simultaneously trained separate networks, the optimization tended to diverge, as illustrated in Figure 3 of \cite{AminiNiaki2021}. To address this issue, we adopt a sequential alternating strategy in which only one network is updated at a time while the remaining networks are kept fixed, and the roles of the networks are then cycled in turn.

For the density $\rho_1$, the loss function is defined as
\begin{equation}
\mathcal{L}_{\rho_1}(\theta_{\rho_1})
=
\mathcal{L}_{\mathrm{ic}}^{\rho_1}
+
\mathcal{L}_{\mathrm{bc}}^{\rho_1}
+
\mathcal{L}_{\mathrm{pde}}^{\rho_1}.
\end{equation}

The initial-condition term is evaluated at initial collocation points $\{x_{\mathrm{ic}}^i\}_{i=1}^{N_{\mathrm{ic}}}$:
\begin{equation}
\mathcal{L}_{\mathrm{ic}}^{\rho_1}
=
\frac{1}{N_{\mathrm{ic}}}
\sum_{i=1}^{N_{\mathrm{ic}}}
\left|
\mathcal{N}_{\rho_1}(x_{\mathrm{ic}}^i,0;\theta_{\rho_1})
-
\rho_{1,0}(x_{\mathrm{ic}}^i)
\right|^2.
\end{equation}

The boundary loss enforces the homogeneous Neumann condition on $\partial\Omega$ using boundary collocation points $\{(x_{\mathrm{bc}}^j,t_{\mathrm{bc}}^j)\}_{j=1}^{N_{\mathrm{bc}}}$:
\begin{equation}
\mathcal{L}_{\mathrm{bc}}^{\rho_1}
=
\frac{1}{N_{\mathrm{bc}}}
\sum_{j=1}^{N_{\mathrm{bc}}}
\left|
\nabla \mathcal{N}_{\rho_1}(x_{\mathrm{bc}}^j,t_{\mathrm{bc}}^j;\theta_{\rho_1})
\cdot \nu(x_{\mathrm{bc}}^j)
\right|^2.
\end{equation}

The PDE loss is evaluated at interior collocation points $\{(x_f^k,t_f^k)\}_{k=1}^{N_f}$:
\begin{equation}
\mathcal{L}_{\mathrm{pde}}^{\rho_1}
=
\frac{1}{N_f}
\sum_{k=1}^{N_f}
\left|
r_{\rho_1}(x_f^k,t_f^k)
\right|^2.
\end{equation}

Similarly, the total losses for $\rho_2$ and $c$ are defined in the same way, with initial-condition, boundary, and PDE residual terms constructed analogously to those of $\rho_1$.

The overall training procedure alternates between the three networks until convergence. At each stage, one network is updated for several optimization steps while the other two remain fixed. We first perform alternating Adam updates and then apply alternating L-BFGS fine-tuning. The procedure is summarized in Algorithm~\ref{alg:pinn_alternating}.

\begin{algorithm}[H]
\SetKwInput{KwInput}{Input}
\SetKwInput{KwOutput}{Output}
\KwInput{Initial data $\rho_{1,0},\rho_{2,0},c_0$, domain $\Omega\times[0,T]$, network architectures, Adam iterations $M$, epochs $K$, learning rate $\eta$, L-BFGS iterations $M_{\text{bfgs}}$}
\KwOutput{Optimized parameters $\theta_{\rho_1}^*,\theta_{\rho_2}^*,\theta_c^*$}
Initialize $\theta_{\rho_1},\theta_{\rho_2},\theta_c$\;
\text{Phase 1: Alternating Adam training}\;
\For{$m=1$ \KwTo $M$}{
    \For{$u \in \{\rho_1,\rho_2,c\}$}{
        Freeze the other two networks\;
        \For{$k=1$ \KwTo $K$}{
            Sample initial, boundary, and interior collocation points\;
            Compute the loss $\mathcal{L}_u$\;
            Update $\theta_u \leftarrow \theta_u-\eta\nabla_{\theta_u}\mathcal{L}_u$\;
        }
    }
}
\text{Phase 2: Alternating L-BFGS fine-tuning}\;
\For{$m=1$ \KwTo $M_{\text{bfgs}}$}{
    \For{$u \in \{\rho_1,\rho_2,c\}$}{
        Freeze the other two networks\;
        Minimize $\mathcal{L}_u(\theta_u)$ using L-BFGS\;
    }
}
\caption{Two-stage alternating training of the PINN model.}
\label{alg:pinn_alternating}
\end{algorithm}

\subsection{A Discrete-Time PINN with Adaptive Collocation for Singular Chemotaxis Dynamics}
\label{sec:discrete_pinn}

In supercritical regimes, solutions of \eqref{eq:full_system_intro} may develop finite-time blow-up as stated in Section~\ref{blow}. While the standard feedforward architecture of Subsection~\ref{sec:continuous_pinn} is effective for many PDE problems and yields reasonable results elsewhere in this work, it becomes difficult to train near blow-up: stiffness, steep gradients, highly localized structures, and competing loss terms demand a more robust formulation. To address this, we adopt a discrete-time PINN framework combining backward Euler time stepping, logarithmic transformations, alternating optimization, and adaptive collocation, together with a modified network architecture.

Specifically, we employ a modified multilayer perceptron (MLP) with Fourier feature input encoding,
\[
\gamma(\mathbf{x}) = [\cos(2\pi \mathbf{B}\mathbf{x}),\, \sin(2\pi \mathbf{B}\mathbf{x})]^\top,
\]
where \(\mathbf{B}\) is a matrix of randomly sampled frequencies. This lifts the input into a higher-dimensional sinusoidal space, letting the network represent high-frequency, sharply localized functions more effectively \cite{Tancik2020}. Two encoder branches further gate every hidden layer via a learned residual/skip connection, adaptively weighting the contribution of each layer \cite{Wang2021}. This gating reduces the gradient pathologies that arise when balancing competing loss terms, such as the PDE residual against initial/boundary conditions, improving training stability near singular structures.

The advantage of the proposed MLP architecture over the standard feedforward tanh architecture within the discrete PINN framework is evident in the sensitivity study of Example~D (Table~\ref{tab:alpha_sensitivity}). Both architectures are identical in all other aspects, they share the same logarithmic transformation, backward Euler time stepping, adaptive collocation strategy, loss functions, and training procedure. The only difference is the network architecture itself: the standard feedforward tanh architecture versus the proposed MLP architecture with Fourier feature encoding and gated residual connections. The modified MLP consistently achieves higher peak values across all values of the clustering exponent \(\alpha\), confirming that the Fourier feature encoding and gated residual connections are beneficial for resolving the highly localized structures characteristic of blow-up, justifying their use in the discrete PINN framework.

We now present the discretization and the corresponding PINN formulation for each case.

\medskip

\noindent
\textbf{Case 1: Parabolic-parabolic regime ($\varepsilon=1$).} 
For this case, we introduce the log-transformed variables for all three unknowns:
\[
\tilde{\rho}_1=\ln(\rho_1+\delta), \qquad
\tilde{\rho}_2=\ln(\rho_2+\delta), \qquad
\tilde{c}=\ln(c+\delta),
\]
with $\delta>0$ small. This transformation ensures positivity, compresses the dynamic range, and improves numerical stability in regimes with strong gradients.

We consider the space-time domain $\Omega_T := \Omega \times [0,T]$ and introduce a uniform partition of the time interval
\[
0 = t^0 < t^1 < \cdots < t^N = T, 
\qquad \Delta t = t^n - t^{n-1}.
\]
For each time level $t^n$, we denote
\[
\rho_i^n(x) \approx \rho_i(x,t^n)~~\text{for}~i=1,2, 
~~\text{and}~~
c^n(x) \approx c(x,t^n)~~
\text{for}~~ x \in \Omega.
\]

Applying backward Euler discretization with step $\Delta t$ to \eqref{eq:full_system_intro} for $\varepsilon=1$ and using the log-transformed variables, we obtain for $n\ge1$:
\begin{equation}
\begin{cases}
\displaystyle
\frac{\tilde{\rho}_1^{n}-\tilde{\rho}_1^{n-1}}{\Delta t}
= \mu_1\Big(\Delta \tilde{\rho}_1^{n}+|\nabla \tilde{\rho}_1^{n}|^2\Big)
-\chi_1 e^{-\tilde{\rho}_1^{n}}
\nabla\cdot\Big((e^{\tilde{\rho}_1^{n}}-\delta)e^{\tilde{c}^{n}}\nabla \tilde{c}^{n}\Big),
\\[6pt]
\displaystyle
\frac{\tilde{\rho}_2^{n}-\tilde{\rho}_2^{n-1}}{\Delta t}
= \mu_2\Big(\Delta \tilde{\rho}_2^{n}+|\nabla \tilde{\rho}_2^{n}|^2\Big)
-\chi_2 e^{-\tilde{\rho}_2^{n}}
\nabla\cdot\Big((e^{\tilde{\rho}_2^{n}}-\delta)e^{\tilde{c}^{n}}\nabla \tilde{c}^{n}\Big),
\\[6pt]
\displaystyle
\frac{\tilde{c}^{n}-\tilde{c}^{n-1}}{\Delta t}
= D\Big(\Delta \tilde{c}^{n}+|\nabla \tilde{c}^{n}|^2\Big)
+\mathrm{e}^{-\tilde{c}^{n}}
\Big[\alpha_1(e^{\tilde{\rho}_1^{n}}-\delta)
+\alpha_2(e^{\tilde{\rho}_2^{n}}-\delta)
-\beta(e^{\tilde{c}^{n}}-\delta)\Big].
\end{cases}
\label{eq:parabolic_discrete}
\end{equation}
The system is completed with initial conditions
\[
\tilde{\rho}_1^{0}(x)=\ln(\rho_1^0(x)+\delta), \qquad
\tilde{\rho}_2^{0}(x)=\ln(\rho_2^0(x)+\delta), \qquad
\tilde{c}^{0}(x)=\ln(c^0(x)+\delta),
\quad x\in\Omega,
\]
and homogeneous Neumann boundary conditions
\[
\frac{\partial \tilde{\rho}_1^{n}}{\partial \nu}
=
\frac{\partial \tilde{\rho}_2^{n}}{\partial \nu}
=
\frac{\partial \tilde{c}^{n}}{\partial \nu}
=0,
\quad \text{on } \partial\Omega.
\]

At each time step $n$, we approximate the unknowns using three independent neural networks
\[
\mathcal{N}_{\tilde{\rho}_1}(x;\theta_{\tilde{\rho}_1}^n), \qquad
\mathcal{N}_{\tilde{\rho}_2}(x;\theta_{\tilde{\rho}_2}^n), \qquad
\mathcal{N}_{\tilde{c}}(x;\theta_{\tilde{c}}^n),
\]
where $\theta_{\tilde{\rho}_1}^n$, $\theta_{\tilde{\rho}_2}^n$, and $\theta_{\tilde{c}}^n$ denote the corresponding sets of weights and biases. The original physical variables are reconstructed through the exponential mapping
\[
\rho_i^n(x)=e^{\tilde{\rho}_i^n(x)}-\eta, \quad i=1,2,
\qquad
c^n(x)=e^{\tilde{c}^n(x)}-\eta,
\]
where $\eta>0$ is small enough to ensure the needed positivity.

The discrete residuals for the parabolic-parabolic case are defined by substituting the neural network approximations into the above scheme, such that
\begin{align}
r_{\tilde{\rho}_1}^{n}(x;\theta_{\tilde{\rho}_1}^{n},\theta_{\tilde{c}}^{n})
&=
\frac{\mathcal{N}_{\tilde{\rho}_1}^{n}(x) - \tilde{\rho}_1^{n-1}(x)}{\Delta t}
-\mu_1\Big(\Delta \mathcal{N}_{\tilde{\rho}_1}^{n}(x) + |\nabla \mathcal{N}_{\tilde{\rho}_1}^{n}(x)|^2\Big) \notag \\
&\quad +\chi_1 e^{-\mathcal{N}_{\tilde{\rho}_1}^{n}(x)}
\nabla\cdot\Big((e^{\mathcal{N}_{\tilde{\rho}_1}^{n}(x)} - \delta)
e^{\mathcal{N}_{\tilde{c}}^{n}(x)}
\nabla \mathcal{N}_{\tilde{c}}^{n}(x)\Big).
\label{eq:residual_rho1_parabolic}
\end{align}
Analogous expressions define $r_{\tilde{\rho}_2}^{n}$ and $r_{\tilde{c}}^{n}$:
\begin{align}
r_{\tilde{\rho}_2}^{n}(x;\theta_{\tilde{\rho}_2}^{n},\theta_{\tilde{c}}^{n})
&=
\frac{\mathcal{N}_{\tilde{\rho}_2}^{n}(x) - \tilde{\rho}_2^{n-1}(x)}{\Delta t}
-\mu_2\Big(\Delta \mathcal{N}_{\tilde{\rho}_2}^{n}(x) + |\nabla \mathcal{N}_{\tilde{\rho}_2}^{n}(x)|^2\Big) \notag \\
&\quad +\chi_2 e^{-\mathcal{N}_{\tilde{\rho}_2}^{n}(x)}
\nabla\cdot\Big((e^{\mathcal{N}_{\tilde{\rho}_2}^{n}(x)} - \delta)
e^{\mathcal{N}_{\tilde{c}}^{n}(x)}
\nabla \mathcal{N}_{\tilde{c}}^{n}(x)\Big),
\label{eq:residual_rho2_parabolic}
\end{align}
and
\begin{align}
r_{\tilde{c}}^{n}(x;\theta_{\tilde{c}}^{n},\theta_{\tilde{\rho}_1}^{n},\theta_{\tilde{\rho}_2}^{n})
&=
\frac{\mathcal{N}_{\tilde{c}}^{n}(x) - \tilde{c}^{n-1}(x)}{\Delta t}
-D\Big(\Delta \mathcal{N}_{\tilde{c}}^{n}(x) + |\nabla \mathcal{N}_{\tilde{c}}^{n}(x)|^2\Big) \notag \\
&\quad - e^{-\mathcal{N}_{\tilde{c}}^{n}(x)}
\Big[\alpha_1(e^{\mathcal{N}_{\tilde{\rho}_1}^{n}(x)} - \delta)
+\alpha_2(e^{\mathcal{N}_{\tilde{\rho}_2}^{n}(x)} - \delta)
-\beta(e^{\mathcal{N}_{\tilde{c}}^{n}(x)} - \delta)\Big].
\label{eq:residual_c_parabolic}
\end{align}

The loss functions for this case are
\[
\mathcal{L}_u^{n}(\theta_u^n)
=
\mathcal{L}_{\mathrm{pde}}^{u,n}(\theta_u^n)
+
\mathcal{L}_{\mathrm{bc}}^{u,n}(\theta_u^n),
\qquad
u \in \{\tilde{\rho}_1,\tilde{\rho}_2,\tilde{c}\},
\]
with
\[
\mathcal{L}_{\mathrm{pde}}^{u,n}(\theta_u^n)
=
\frac{1}{N_f}\sum_{k=1}^{N_f}
\left|r_{u}^{n}(x_k;\theta_u^n,\theta_{\mathrm{coupled}}^n)\right|^2,
\]
and
\[
\mathcal{L}_{\mathrm{bc}}^{u,n}(\theta_u^n)
=
\frac{1}{N_{\mathrm{bc}}}\sum_{j=1}^{N_{\mathrm{bc}}}
\left|\nabla \mathcal{N}_{u}(x_{\mathrm{bc}}^j;\theta_u^n)\cdot \nu(x_{\mathrm{bc}}^j)\right|^2.
\]
Here, $\theta_{\mathrm{coupled}}^n$ denotes the parameters of the other networks, which are kept fixed during each alternating update.
\medskip

\noindent
\noindent
\textbf{Case 2: Parabolic-elliptic regime ($\varepsilon=0$).} 
In this limit, the third equation reduces to an elliptic relation. While the density variables retain their logarithmic form, the chemoattractant $c$ is not log-transformed:
\[
\tilde{\rho}_1=\ln(\rho_1+\delta), \qquad
\tilde{\rho}_2=\ln(\rho_2+\delta), \qquad
\tilde{c}=c.
\]
This choice is deliberate. When $\varepsilon=0$, the chemoattractant equation becomes
\[
D\Delta c = \beta c - \alpha_1 \rho_1 - \alpha_2 \rho_2,
\]
which is linear in $c$ with constant coefficients. Applying a logarithmic transformation $\tilde{c} = \ln(c+\delta)$ would turn this into the severely nonlinear expression
\[
D e^{\tilde{c}}(\Delta \tilde{c} + |\nabla \tilde{c}|^2) = \beta e^{\tilde{c}} - \alpha_1(\cdots) - \alpha_2(\cdots),
\]
destroying the well-conditioned linear elliptic structure and introducing exponential nonlinearities that degrade convergence and training stability. Thus, keeping $c$ in its native form is essential for maintaining the favorable properties of the elliptic subsystem. 

In contrast, when \(\varepsilon = 1\), the chemoattractant equation is parabolic and \(c\) evolves in time. Although \(c\) itself does not blow up, it can develop steep gradients in response to the rapidly growing densities, which can degrade numerical stability. Applying the same logarithmic transformation to \(c\) as to the densities enforces positivity, compresses the dynamic range, and improves the conditioning of the parabolic system. This is why we treat \(c\) differently in the two regimes.

Next, applying backward Euler discretization to the density equations and solving the elliptic equation at each time level for $\varepsilon=0$ yields
\begin{equation}
\begin{cases}
\displaystyle
\frac{\tilde{\rho}_1^{n}-\tilde{\rho}_1^{n-1}}{\Delta t}
= \mu_1\Big(\Delta \tilde{\rho}_1^{n}+|\nabla \tilde{\rho}_1^{n}|^2\Big)
-\chi_1 e^{-\tilde{\rho}_1^{n}}
\nabla\cdot\Big((e^{\tilde{\rho}_1^{n}}-\delta)\nabla c^{n}\Big),
\\[6pt]
\displaystyle
\frac{\tilde{\rho}_2^{n}-\tilde{\rho}_2^{n-1}}{\Delta t}
= \mu_2\Big(\Delta \tilde{\rho}_2^{n}+|\nabla \tilde{\rho}_2^{n}|^2\Big)
-\chi_2 e^{-\tilde{\rho}_2^{n}}
\nabla\cdot\Big((e^{\tilde{\rho}_2^{n}}-\delta)\nabla c^{n}\Big),
\\[6pt]
\displaystyle
0 = D\Delta c^{n}
+\alpha_1(e^{\tilde{\rho}_1^{n}}-\delta)
+\alpha_2(e^{\tilde{\rho}_2^{n}}-\delta)
-\beta c^{n}.
\end{cases}
\label{eq:elliptic_discrete}
\end{equation}
The initial and boundary conditions for $\tilde{\rho}_1$ and $\tilde{\rho}_2$ are the same as in the parabolic-parabolic case:
\[
\tilde{\rho}_1^{0}(x)=\ln(\rho_1^0(x)+\delta), \qquad
\tilde{\rho}_2^{0}(x)=\ln(\rho_2^0(x)+\delta), \qquad
\frac{\partial \tilde{\rho}_1^n}{\partial \nu}
=
\frac{\partial \tilde{\rho}_2^n}{\partial \nu}
=0,
\quad \text{on } \partial\Omega,
\]
while for $c$ we impose
\[
c^0(x)=c_0(x), \qquad
\frac{\partial c^n}{\partial \nu}=0 \quad \text{on } \partial\Omega.
\]

At each time step, we approximate the unknowns using three independent neural networks
\[
\mathcal{N}_{\tilde{\rho}_1}(x;\theta_{\tilde{\rho}_1}^n), \qquad
\mathcal{N}_{\tilde{\rho}_2}(x;\theta_{\tilde{\rho}_2}^n), \qquad
\mathcal{N}_{c}(x;\theta_{c}^n).
\]
The original physical densities are reconstructed as
\[
\rho_i^n(x)=e^{\mathcal{N}_{\tilde{\rho}_i}^n(x)}-\eta, \qquad i=1,2,
\]
where $\eta>0$ is a small constant.

The residuals for the parabolic-elliptic case are defined by substituting the neural network approximations into the above scheme:
\begin{align}
r_{\tilde{\rho}_1}^{n,0}(x;\theta_{\tilde{\rho}_1}^n,\theta_c^n)
&=
\frac{\mathcal{N}_{\tilde{\rho}_1}^n(x) - \tilde{\rho}_1^{n-1}(x)}{\Delta t}
-\mu_1\Big(\Delta \mathcal{N}_{\tilde{\rho}_1}^n(x) + |\nabla \mathcal{N}_{\tilde{\rho}_1}^n(x)|^2\Big) \notag \\
&\quad +\chi_1 e^{-\mathcal{N}_{\tilde{\rho}_1}^n(x)}
\nabla\cdot\Big((e^{\mathcal{N}_{\tilde{\rho}_1}^n(x)}-\delta)\nabla \mathcal{N}_c^n(x)\Big),
\label{eq:residual_rho1_elliptic}
\\[4pt]
r_{\tilde{\rho}_2}^{n,0}(x;\theta_{\tilde{\rho}_2}^n,\theta_c^n)
&=
\frac{\mathcal{N}_{\tilde{\rho}_2}^n(x) - \tilde{\rho}_2^{n-1}(x)}{\Delta t}
-\mu_2\Big(\Delta \mathcal{N}_{\tilde{\rho}_2}^n(x) + |\nabla \mathcal{N}_{\tilde{\rho}_2}^n(x)|^2\Big) \notag \\
&\quad +\chi_2 e^{-\mathcal{N}_{\tilde{\rho}_2}^n(x)}
\nabla\cdot\Big((e^{\mathcal{N}_{\tilde{\rho}_2}^n(x)}-\delta)\nabla \mathcal{N}_c^n(x)\Big),
\label{eq:residual_rho2_elliptic}
\intertext{and}
r_c^{n,0}(x;\theta_c^n,\theta_{\tilde{\rho}_1}^n,\theta_{\tilde{\rho}_2}^n)
&=
D\Delta \mathcal{N}_c^n(x)
+\alpha_1(e^{\mathcal{N}_{\tilde{\rho}_1}^n(x)}-\delta)
+\alpha_2(e^{\mathcal{N}_{\tilde{\rho}_2}^n(x)}-\delta)
-\beta \mathcal{N}_c^n(x).
\label{eq:residual_c_elliptic}
\end{align}

The loss functions for $\tilde{\rho}_1$ and $\tilde{\rho}_2$ are defined as in the parabolic-parabolic case:
\[
\mathcal{L}_{\tilde{\rho}_i}^{n,0}(\theta_{\tilde{\rho}_i}^n)
=
\frac{1}{N_f}\sum_{k=1}^{N_f}
\left|r_{\tilde{\rho}_i}^{n,0}(x_k;\theta_{\tilde{\rho}_i}^n,\theta_c^n)\right|^2
+
\mathcal{L}_{\mathrm{bc}}^{\tilde{\rho}_i,n,0}(\theta_{\tilde{\rho}_i}^n),
\qquad i=1,2,
\]
with
\[
\mathcal{L}_{\mathrm{bc}}^{\tilde{\rho}_i,n,0}(\theta_{\tilde{\rho}_i}^n)
=
\frac{1}{N_{\mathrm{bc}}}\sum_{j=1}^{N_{\mathrm{bc}}}
\left|\nabla \mathcal{N}_{\tilde{\rho}_i}(x_{\mathrm{bc}}^j;\theta_{\tilde{\rho}_i}^n)\cdot \nu(x_{\mathrm{bc}}^j)\right|^2.
\]

\begin{figure}[t]
\centering
\includegraphics[width=0.5\textwidth]{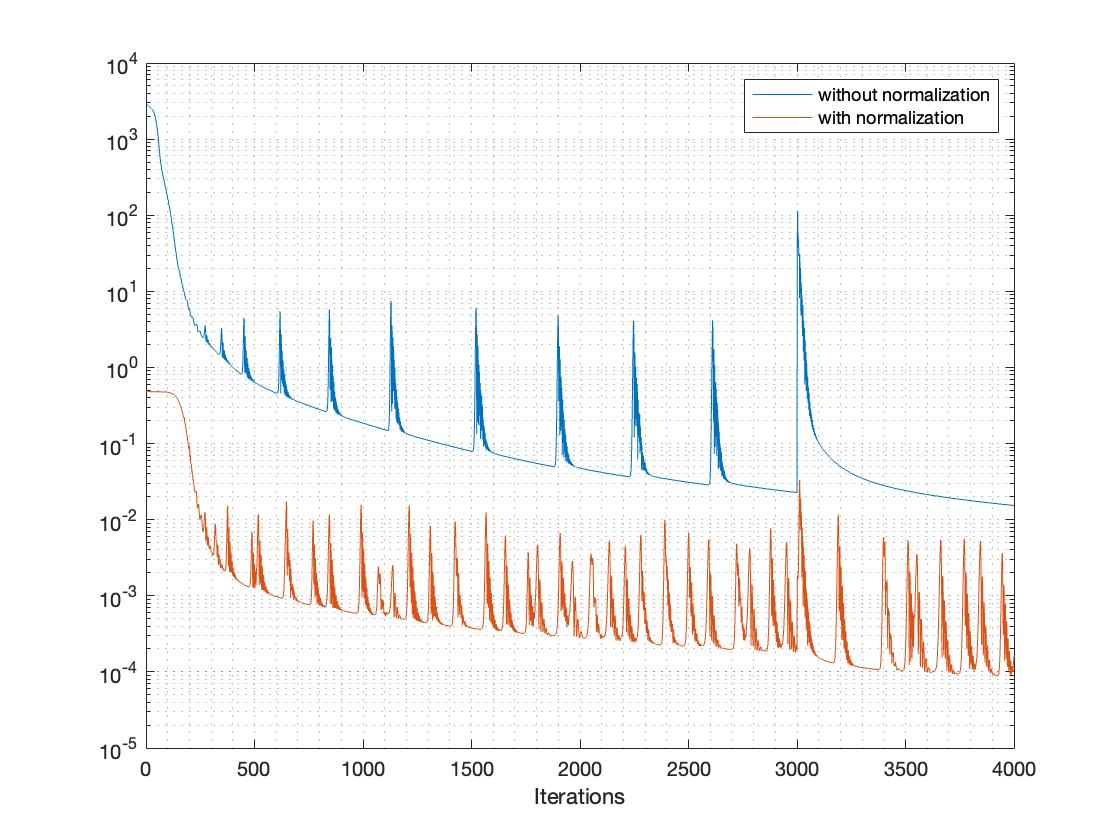}
\caption{Evolution of the chemoattractant residual loss $\mathcal{L}_c^{n,0}$ for the elliptic equation with and without source-scale residual normalization.}
\label{fig:loss_comparison}
\end{figure}

For $c$, however, a standard MSE loss would be dominated by regions where the source terms $\alpha_1\rho_1+\alpha_2\rho_2$ are extremely large, as these can span several orders of magnitude. To prevent this imbalance, we normalize the residual by the local source scale. The loss for $c$ is therefore defined as
\[
\mathcal{L}_{c}^{n,0}(\theta_c^n)
=
\frac{1}{N_f}\sum_{k=1}^{N_f}
\left(
\frac{
r_c^{n,0}(x_k;\theta_c^n,\theta_{\tilde{\rho}_1}^n,\theta_{\tilde{\rho}_2}^n)
}{
\exp\bigl(\mathcal{N}_{\tilde{\rho}_1}(x_k;\theta_{\tilde{\rho}_1}^n)\bigr)
+
\exp\bigl(\mathcal{N}_{\tilde{\rho}_2}(x_k;\theta_{\tilde{\rho}_2}^n)\bigr)
+ 1
}
\right)^2
+
\mathcal{L}_{\mathrm{bc}}^{c,n,0}(\theta_c^n),
\]
where the boundary loss for $c$ is given by
\[
\mathcal{L}_{\mathrm{bc}}^{c,n,0}(\theta_c^n)
=
\frac{1}{N_{\mathrm{bc}}}\sum_{j=1}^{N_{\mathrm{bc}}}
\left|\nabla \mathcal{N}_{c}(x_{\mathrm{bc}}^j;\theta_c^n)\cdot \nu(x_{\mathrm{bc}}^j)\right|^2.
\]
The normalization by $\exp(\tilde{\rho}_1) + \exp(\tilde{\rho}_2)$ balances the contribution of each collocation point to the total loss, ensuring that the network resolves the elliptic equation uniformly across the domain while still capturing the sharp variations induced by the densities. The $+1$ term provides numerical stability when both densities are negligible. This normalized residual formulation is strongly preferred, as it decouples the loss magnitude from the local solution scale, stabilizes training, and yields a more accurate approximation of $c$ in the parabolic-elliptic setting. 
Figure~\ref{fig:loss_comparison} demonstrates this improvement by comparing the evolution of the chemoattractant loss $\mathcal{L}_c^{n,0}$ with and without the source-scale normalization. The normalized loss achieves smoother convergence and significantly lower final values, confirming that pointwise weighting prevents the loss from being dominated by high-magnitude source regions and improves the overall approximation accuracy.

\medskip

\begin{figure}[t]
\centering
\includegraphics[width=0.45\textwidth]{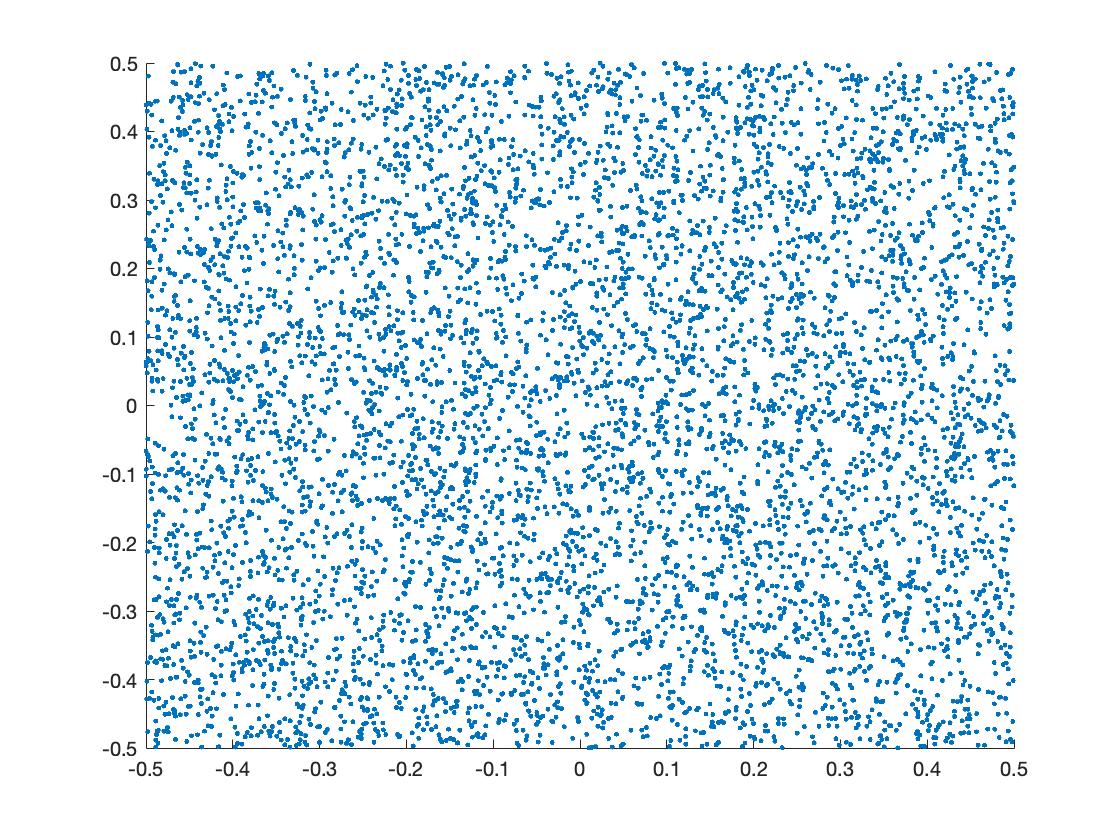}
\hfill
\includegraphics[width=0.45\textwidth]{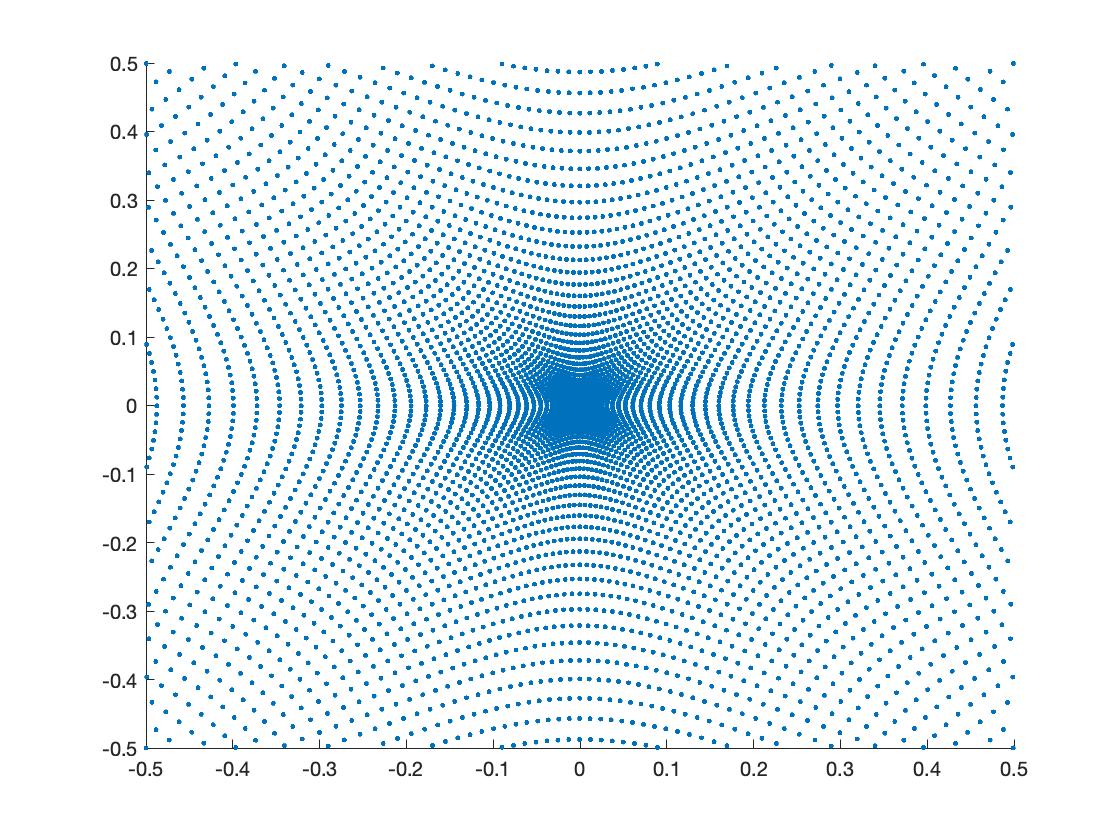}
\caption{Illustration of the adaptive collocation strategy. 
Left: initial uniform candidate points over $\Omega$. 
Right: clustered collocation points after the transformation and recentering at $x_{u,\max}^{n-1}$.}
\label{fig:collocation_cluster}
\end{figure}

To construct an efficient adaptive collocation set for training the neural network approximations at time level $t^n$, we exploit the solution at the previous time level $t^{n-1}$. Since the current solution $u^n$ is unknown and must be learned by the network, we use $u^{n-1}$ as a proxy to identify regions where sharp gradients or singularities are likely to develop. For each variable $u \in \{\tilde{\rho}_1,\tilde{\rho}_2,\tilde{c}\}$ (or $c$ in the elliptic case), we first generate a very fine uniform set of candidate collocation points over $\Omega$, which we denote as $\Omega_{\text{fine}}$. Among these points, we locate the spatial position where the absolute value of the previous-time solution $u^{n-1}(x)$ attains its maximum:
\[
x_{u,\max}^{n-1} = \arg\max_{x \in \Omega_{\text{fine}}} |u^{n-1}(x)|.
\]
This point is then used as the center of a locally refined cluster for training the network corresponding to $u^n$ at the current time step. The cluster is constructed by first sampling points uniformly inside a disk of radius $r$ centered at the origin, and then applying a power-law transformation $(x,y)\mapsto (x^\alpha,y^\alpha)$, with $\alpha>1$, so that the points become increasingly concentrated toward the center. 
A systematic sensitivity study of the clustering exponent $\alpha$ and the selection criterion used in this work are presented in Section~4.2.1 (Example~D). After shifting the resulting cloud to $x_{u,\max}^{n-1}$, we restrict the points to the physical domain $\Omega$ by a filtering step; see Figure~\ref{fig:collocation_cluster} for the collocation points before and after clustering. This procedure yields a smoothly varying non-uniform distribution with density increasing gradually toward the center, effectively placing the highest concentration of collocation points near the maximum of $u^{n-1}$, where blow-up or sharp gradients are expected.

\begin{rem}
Accurate resolution of finite-time blow-up strongly depends on the spatial distribution of collocation points. Residual-based adaptive refinement methods
\cite{Wu2023,MaoMeng2023,Gao2023}
typically place additional points in regions where the PDE residual is largest. In chemotaxis models approaching blow-up, these regions coincide with locations where the solution develops extremely sharp gradients and highly localized peaks. Although such refinement increases the local point density, the resulting point cloud may become irregular and unevenly distributed, making it difficult for the neural network to accurately represent the rapidly varying solution profile.

To address this issue, we employ a smooth and symmetric clustering strategy centered around the predicted blow-up location. The point density increases gradually toward the region of highest stiffness and decreases smoothly away from it, avoiding abrupt changes in spatial resolution. This provides a more structured representation of the singular region while maintaining adequate coverage of the remainder of the domain. As illustrated in Figure 3, both approaches are compared at the same final time. It is important to note that both the residual-based adaptive refinement strategy and the proposed smooth clustering strategy are configured with the same total number of collocation points and identical neural network parameters, ensuring a fair and direct comparison. The residual-based refinement underestimates the peak magnitude, whereas the proposed smooth clustering accurately captures the highly concentrated solution profile.
\end{rem}

\begin{figure}[t]
\centering
\begin{minipage}{0.45\textwidth}
    \centering
    \includegraphics[width=\linewidth]{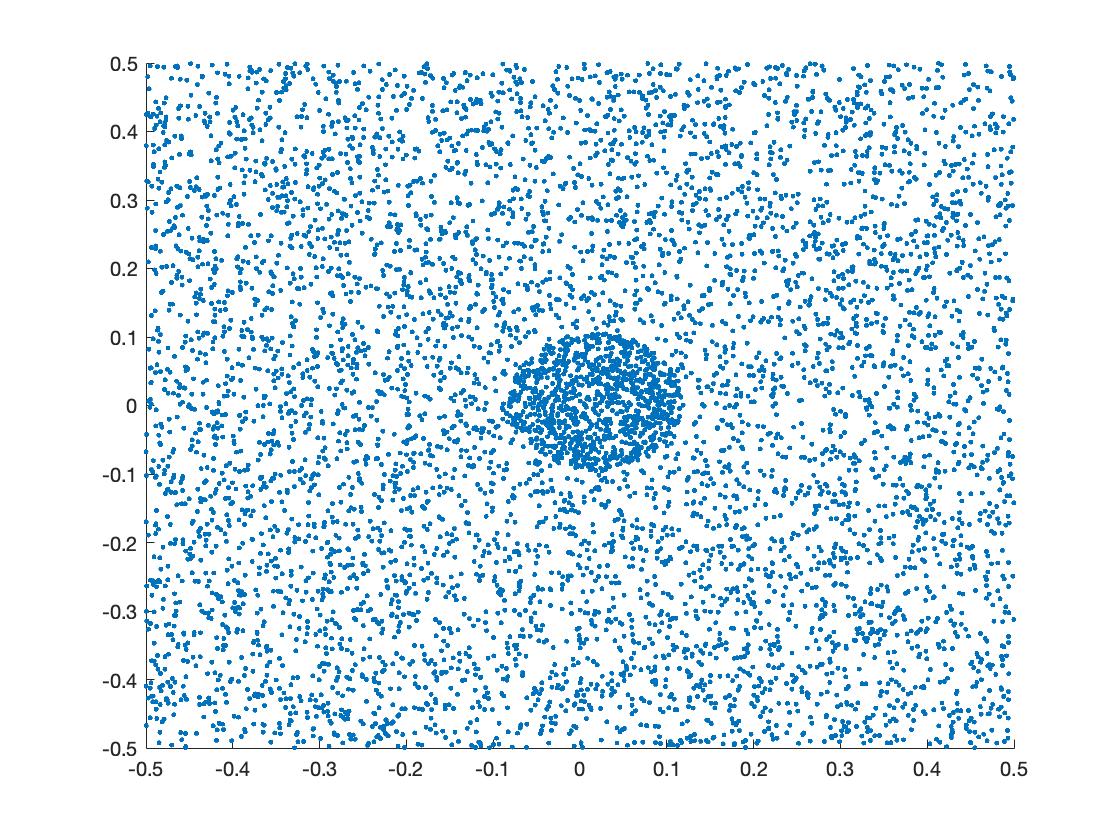}
    \caption*{(a) Residual-based adaptive points}
\end{minipage}
\hfill
\begin{minipage}{0.45\textwidth}
    \centering
    \includegraphics[width=\linewidth]{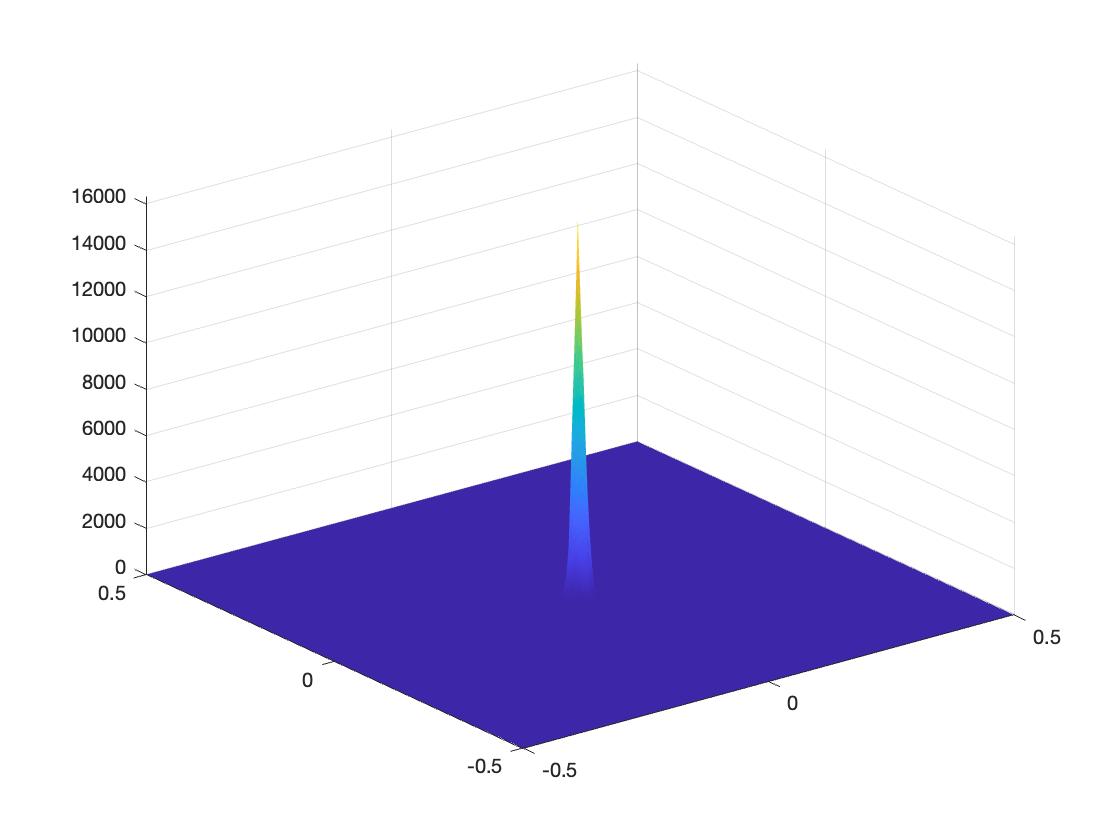}
    \caption*{(b) Discrete PINN solution at final time}
\end{minipage}

\vspace{0.4cm}

\begin{minipage}{0.45\textwidth}
    \centering
    \includegraphics[width=\linewidth]{figures/example_d_smooth/smooth_points.jpg}
    \caption*{(c) Smooth clustered points (proposed)}
\end{minipage}
\hfill
\begin{minipage}{0.45\textwidth}
    \centering
    \includegraphics[width=\linewidth]{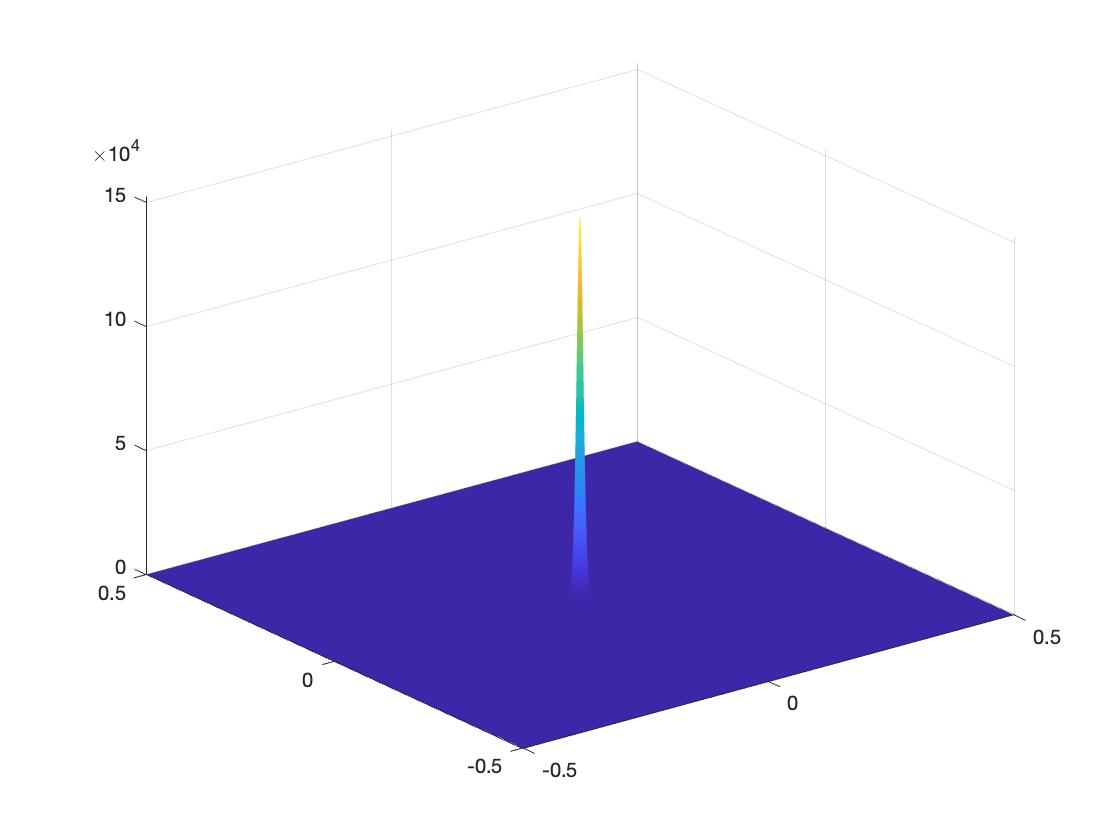}
    \caption*{(d) Discrete PINN solution at same final time}
\end{minipage}

\caption{Effect of collocation distribution on blow-up approximation. Smooth clustering resolves the peak more sharply and yields a higher maximum value than residual-based adaptive refinement at the same final time.}
\label{fig:smooth_vs_irregular}
\end{figure}

After determining the adaptive set of collocation points for the current time step, we proceed with the training of the neural networks. At each time step, the optimization is performed using an alternating strategy: one network is updated while the other two are kept fixed. The training is carried out in two phases. First, we apply AdamW \cite{loshchilov2019adamw} for global exploration of the loss landscape. As an adaptive stochastic gradient method, AdamW adjusts the learning rate of each parameter individually, assigning smaller effective updates to parameters associated with consistently large gradients, which helps stabilize training and accelerates convergence in the early stages. Once the loss reaches a sufficiently low plateau, we switch to L-BFGS, a quasi-Newton method, to perform accurate local refinement. This combination of global exploration (AdamW) and local optimization (L-BFGS) is particularly effective in our setting, as it allows the networks to capture the sharp gradients and potential blow-up behavior induced by the clustered collocation points without getting trapped in poor local minima.

To ensure numerical stability in the presence of rapidly growing solutions, we implement a conservative rollback strategy. At the beginning of each time step, a full checkpoint of all network parameters is saved. If any loss value or gradient becomes non-finite during the Adam or L-BFGS optimization phases, the training is rolled back to the last valid state. This prevents catastrophic divergence and allows the simulation to continue even when the solution approaches the blow-up threshold. This rollback mechanism is essential for the robustness of the method near singularities. The overall procedure, including the clustering strategy and the alternating optimization, is summarized in Algorithm~2.

\begin{algorithm}[t]
\small
\SetKwInput{KwInput}{Input}
\SetKwInput{KwOutput}{Output}

\KwInput{Initial data $\tilde{\rho}_1^0,\tilde{\rho}_2^0,\tilde{c}^0$, fine domain $\Omega_{\text{fine}}$, number of time steps $N$, AdamW alternations $K_{\text{alt}}$, L-BFGS alternations $K_{\text{bfgs}}$}

\KwOutput{Networks $\theta_{\tilde{\rho}_1}^N,\theta_{\tilde{\rho}_2}^N,\theta_{\tilde{c}}^N$}

Initialize $\tilde{\rho}_1^0(x), \tilde{\rho}_2^0(x), \tilde{c}^0(x)$ from the initial conditions\;

\For{$n=1$ \KwTo $N$}{

    \textbf{Adaptive collocation construction}\;

    \For{$u \in \{\tilde{\rho}_1,\tilde{\rho}_2,\tilde{c}\}$}{

        Compute
        $x_{u,\max}^{n-1}
        = \arg\max_{x\in\Omega_{\text{fine}}}|u^{n-1}(x)|$\;

        Generate clustered collocation points
        $\mathcal{C}_u^n$
        around $x_{u,\max}^{n-1}$ using the power-law mapping\;
    }

    Initialize
    $\mathcal{N}_{\tilde{\rho}_1}^n$,
    $\mathcal{N}_{\tilde{\rho}_2}^n$,
    $\mathcal{N}_{\tilde{c}}^n$\;

    \textbf{Phase 1: Alternating AdamW training}\;

    \For{$k=1$ \KwTo $K_{\text{alt}}$}{

        \For{$u \in \{\tilde{\rho}_1,\tilde{\rho}_2,\tilde{c}\}$}{

            Freeze the other two networks\;

            Compute $\mathcal{L}_u^n$\;

            Update $\theta_u^n$ via AdamW on
            $\mathcal{C}_u^n \cup \mathcal{C}_{\mathrm{bc}}$\;
        }
    }

    \textbf{Phase 2: Alternating L-BFGS fine-tuning}\;

    \For{$k=1$ \KwTo $K_{\text{bfgs}}$}{

        \For{$u \in \{\tilde{\rho}_1,\tilde{\rho}_2,\tilde{c}\}$}{

            Freeze the other two networks\;

            Minimize $\mathcal{L}_u^n(\theta_u)$ using L-BFGS\;
        }
    }
}

\caption{Discrete-time PINN with adaptive clustering and alternating AdamW/L-BFGS training. A checkpoint-restore mechanism is used to roll back the optimization whenever non-finite values are encountered.}
\label{alg:discrete_pinn}
\end{algorithm}

\section{Numerical examples}

We evaluate the proposed PINN frameworks on six benchmark problems covering subcritical regimes with smooth dynamics and supercritical regimes with finite-time blow-up. For all examples, we set \(N_f = 8000\) interior collocation points. Hyperparameters were selected using Optuna~\cite{Akiba2019Optuna}; network configurations are in Table~\ref{tab:example_parameters}. Models are trained on a Tesla P100 12GB GPU using Python 3.12.4, with Ray enabling up to ten concurrent runs per GPU. Code is available at \url{https://github.com/ELbahja88/hybrid-pinn-chemotaxis}.

\begin{table}[htbp]
\centering
\caption{Network hyperparameter settings for all numerical examples.}
\label{tab:example_parameters}
\begin{tabular}{c|c|c|c}
\hline
Example & Regime & Network depth & Width \\
\hline
A & Subcritical (1D)                 & 7 & 27  \\
B & Subcritical (2D)                 & 5 & 32  \\
C & Subcritical (full system)        & 6 & 24  \\
D & Supercritical (single-species)   & 4 & 30  \\
E & Supercritical (two-species PE)   & 4 & 30  \\
F & Supercritical (two-species PP)   & 4 & 30  \\
\hline
\end{tabular}
\end{table}

\subsection{Subcritical regime}

We assess the proposed PINN in the subcritical regime, where solutions remain smooth and globally bounded. For these tests, we employ the continuous-time PINN described in Section~3.1, which is computationally efficient and well suited to regular solutions. Validation is performed using manufactured solutions of the chemotaxis system \eqref{eq:full_system_intro} under homogeneous Neumann boundary conditions. Three benchmark cases are considered: a one-dimensional test, a two-dimensional single-species diagnostic problem, and the full system. In all experiments, the physical parameters are set to unity,
\[
\mu_1 = \mu_2 = \chi_1 = \chi_2 = D = \alpha_1 = \alpha_2 = \beta = \varepsilon = 1,
\]
and the source terms are obtained by exact substitution of the analytical solutions into \eqref{eq:full_system_intro}, ensuring zero residual.
\subsubsection{Example A: 1D problem}

We first consider the one-dimensional domain \(x\in[0,\pi]\), \(t\in[0,1]\), with exact solution
\[
\rho_1(x,t)=e^{-t}\cos x,\qquad 
c(x,t)=e^{-t}\cos x,\qquad 
\rho_2 \equiv 0.
\]
The corresponding forcing term is
\[
r_1(x,t) = -e^{-2t}\cos(2x),\qquad r_2(x,t)=0.
\]

Since \(\rho_1=c\), only \(\rho_1\) is reported. Table~\ref{tab:1d_errors} compares the relative $L^2$ and $L^\infty$ errors of the PINN with those of the DDGIC scheme with positivity-preserving limiter on a mesh of \(N=40\) cells using \(P^2\) elements presented in \cite{QiuLiuYan2021}.

\begin{table}[htbp]
\centering
\caption{1D results at \(T=0.2\) (relative errors).}
\label{tab:1d_errors}
\begin{tabular}{c|cc|cc}
\hline
 & \multicolumn{2}{c|}{\cite{QiuLiuYan2021}} & \multicolumn{2}{c}{PINN } \\
Variable & \(L^2\) & \(L^\infty\) & \(L^2\) & \(L^\infty\) \\
\hline
$\rho_1$ & \(2.02\times10^{-3}\) & \(3.96\times10^{-3}\) & \(2.3\times10^{-4}\) & \(4.9\times10^{-4}\) \\
\hline
\end{tabular}
\end{table}

The proposed PINN enforces positivity via a softplus output layer, avoiding any post-processing limiter. This leads to improved accuracy compared with the DG reference solution while preserving the physical constraint intrinsically.

\subsubsection{Example B: 2D single-species problem}

We next consider the two-dimensional domain $\Omega=[0,\pi]^2$ with solution
\[
\rho_1(x,y,t)=e^{-t}(\cos x+\cos y),\qquad 
c(x,y,t)=e^{-t}(\cos x+\cos y),\qquad \rho_2\equiv 0.
\]
The forcing term reads
\[
r_1(x,y,t) = -2e^{-2t}\bigl(\cos^2 x + \cos x\cos y + \cos^2 y - 1\bigr),\qquad r_2\equiv 0.
\]

\begin{table}[htbp]
\centering
\caption{2D single-species results at $T=0.2$ (absolute errors).}
\label{tab:2d_errors}
\begin{tabular}{c|cc|cc}
\hline
 & \multicolumn{2}{c|}{\cite{SulmanNguyen2019}} & \multicolumn{2}{c}{PINN } \\
Variable & $L^2$ error & $L^1$ error & $L^2$ error & $L^1$ error \\
\hline
$\rho_1$ & $1.02\times10^{-1}$ & $4.52\times10^{-1}$ & ${1.48\times10^{-1}}$ & ${1.98\times10^{-1}}$ \\
\hline
\end{tabular}
\end{table}

Table~\ref{tab:2d_errors} compares the absolute \(L^1\) and \(L^2\) errors of the proposed PINN at \(T=0.2\) with those of the moving-mesh finite element method (PMA-FEM) of~\cite{SulmanNguyen2019} on a coarse \(11\times11\) mesh. The PINN achieves competitive accuracy on this mesh, with a slightly smaller \(L^1\) error and a comparable \(L^2\) error, while offering several practical advantages: it inherently preserves positivity without limiting or post-processing, operates mesh-free, and avoids the complexity of mesh movement and reconnection. Although the PMA-FEM can achieve higher accuracy on finer grids, the PINN provides a compelling alternative for problems where rapid prototyping or coarse-resolution simulations suffice, or where mesh generation is challenging.

\subsubsection{Example C: Multi-species problem}

We finally consider the full two-species system \eqref{eq:full_system_intro} on \(\Omega=[0,\pi]^2\). To further assess the capability of the proposed PINN in a fully coupled setting, we introduce the following manufactured solution:
\[
\begin{aligned}
\rho_1(x,y,t) &= e^{-t}\sin x \sin y,\\
\rho_2(x,y,t) &= e^{-t}\cos x \cos y,\\
c(x,y,t) &= e^{-t}\bigl(\sin x \sin y + \cos x \cos y\bigr).
\end{aligned}
\]

Substituting this solution into \eqref{eq:full_system_intro} yields the corresponding source terms:
\[
\begin{aligned}
r_1(x,y,t) &= e^{-t}\sin x\sin y \\
&\quad + e^{-2t}\Bigl[\sin^2 x + \sin^2 y - 4\sin^2 x\sin^2 y - 4\sin x\cos x\sin y\cos y\Bigr],\\[4pt]
r_2(x,y,t) &= e^{-t}\cos x\cos y \\
&\quad + e^{-2t}\Bigl[-4\sin x\cos x\sin y\cos y - \cos(2x)\cos^2 y - \cos^2 x\cos(2y)\Bigr],\\[4pt]
r_3(x,y,t) &= e^{-t}\bigl(\sin x\sin y + \cos x\cos y\bigr).
\end{aligned}
\]

The manufactured solution satisfies the PDE exactly with these source terms. For this verification test, Dirichlet boundary conditions are imposed instead of the original homogeneous Neumann boundary conditions. The experiment is designed to evaluate the capability of PINNs to accurately capture the dynamics of multi-species chemotaxis equations.
\begin{table}[htbp]
\centering
\caption{Relative errors for the multi-species system at final time \(t=1\).}
\label{tab:multi_species_errors}
\begin{tabular}{c|cc}
\hline
Variable &   \(L^2\) error &   \(L^\infty\) error \\
\hline
$\rho_1$ & $1.3\times 10^{-3}$ & $2.4\times10^{-3}$ \\
$\rho_2$ & $9.3\times10^{-4}$ & $2.7\times10^{-3}$ \\
$c$      & $2\times10^{-3}$ & $4.5\times10^{-3}$ \\
\hline
\end{tabular}
\end{table}

Table~\ref{tab:multi_species_errors} shows that the relative errors are consistently small for all variables, with both relative $L^2$ and $L^\infty$ errors of order $10^{-3}$. This indicates that the PINN provides an accurate approximation of the multi-species chemotaxis system and effectively captures the coupling between the different species and the chemoattractant.

\subsection{Supercritical regimes}

We now apply the discrete PINN framework introduced in Section~\ref{sec:discrete_pinn} to supercritical regimes, where the original system~\eqref{eq:full_system_intro} is known to develop finite-time singularities. In all examples, we employ a logarithmic transformation with \(\delta = 10^{-6}\), together with the residual-adaptive clustering strategy and the backward Euler time discretization described in Algorithm~\ref{alg:discrete_pinn}. The hyperparameters are chosen according to Table~\ref{tab:example_parameters}. We use \(\alpha = 6\) for all supercritical examples, with \(r = 1\) for Example D and \(r = 3\) for Examples E and F.

\subsubsection{Example D: Single-species Keller--Segel blow-up}

We consider the single-species regime of~\eqref{eq:full_system_intro}, obtained by setting 
\(\rho_2 \equiv 0\) and \(\varepsilon = 1\) (fully parabolic case). The parameters are
\[
\mu_1 = 1,\quad \chi_1 = 1,\quad D = 1,\quad \alpha_1 = 1,\quad \beta = 1.
\]
The spatial domain is \(\Omega = [-0.5,0.5]^2\), equipped with homogeneous Neumann boundary conditions.

The initial data are given by centrally concentrated Gaussians:
\begin{align}
\rho_1(x,y,0) &= 840\,e^{-84(x^2+y^2)},\\
c(x,y,0) &= 420\,e^{-42(x^2+y^2)}, \qquad (x,y)\in\Omega.
\end{align}
The total mass satisfies \(M = 10\pi\), which exceeds the critical threshold \(8\pi\). 
According to the results recalled in Section~\ref{blow}, radially symmetric solutions 
in this regime are expected to develop finite-time blow-up at the origin.

\begin{figure}[t]
\centering
\begin{minipage}{0.48\textwidth}
    \centering
    \includegraphics[width=\linewidth]{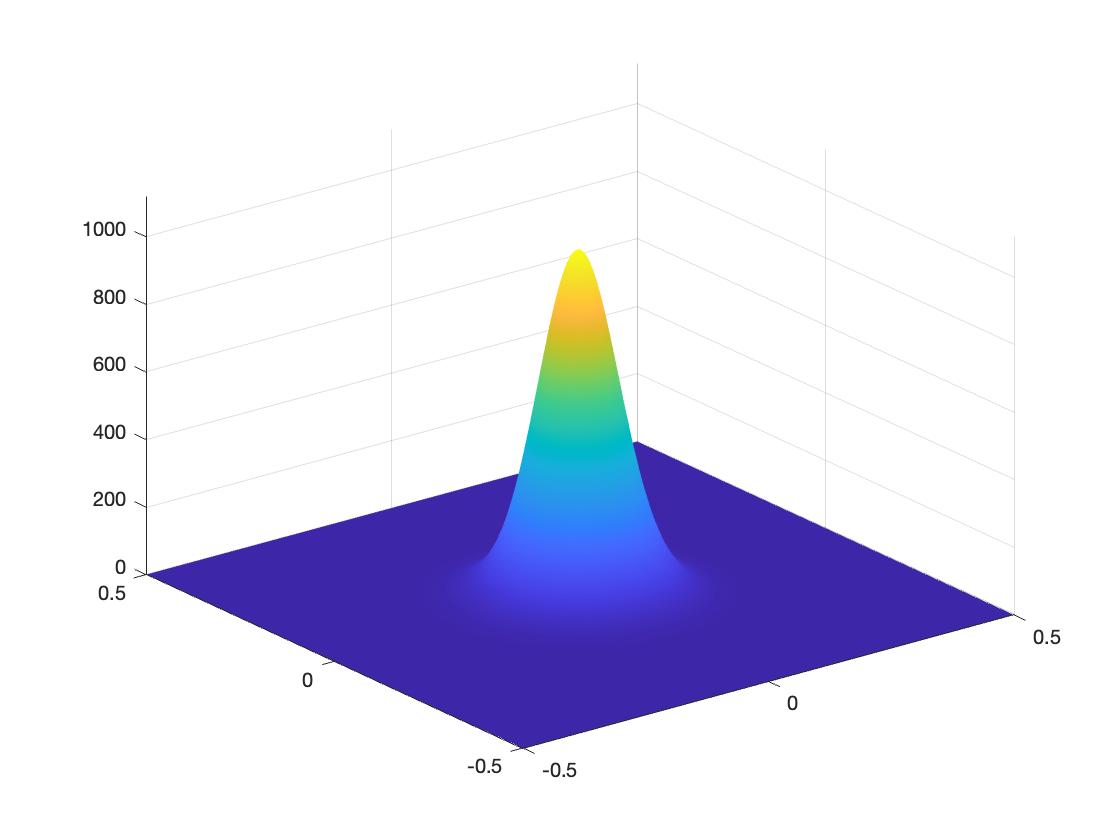}
    \caption*{$t = 3\times 10^{-6}$}
\end{minipage}
\hfill
\begin{minipage}{0.48\textwidth}
    \centering
    \includegraphics[width=\linewidth]{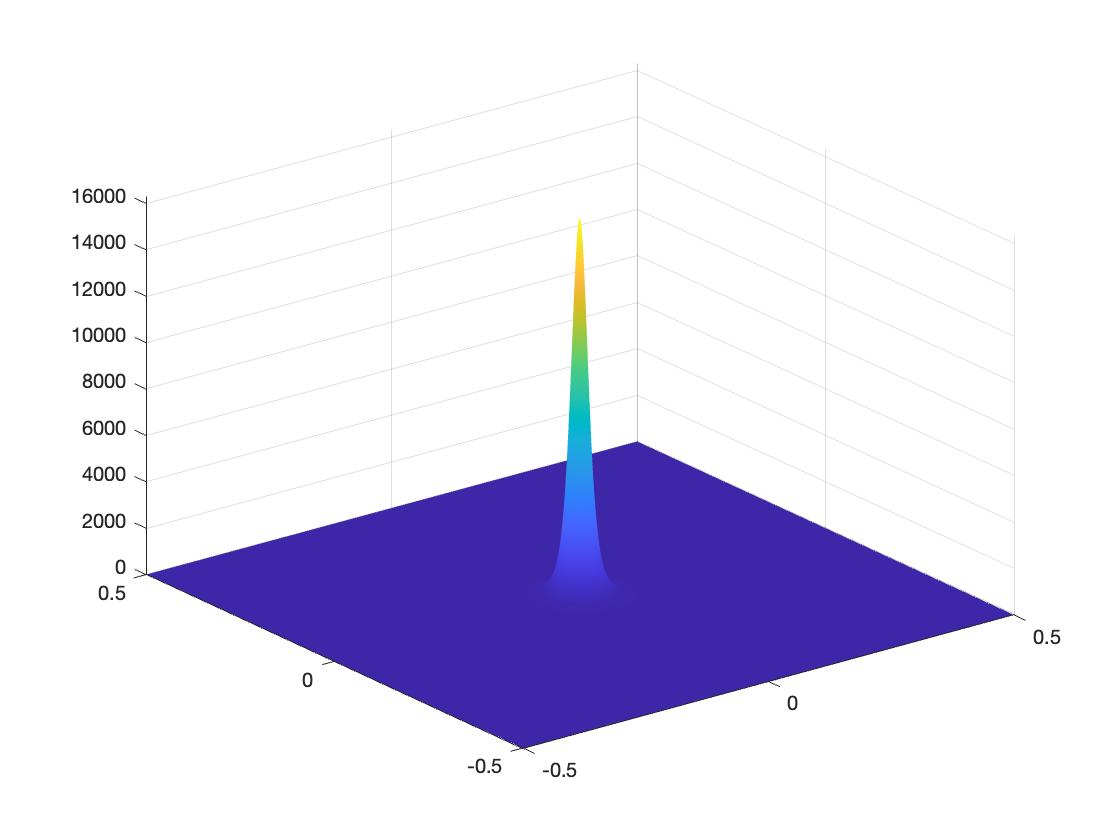}
    \caption*{$t = 5\times 10^{-5}$}
\end{minipage}

\vspace{0.5em}

\begin{minipage}{0.48\textwidth}
    \centering
    \includegraphics[width=\linewidth]{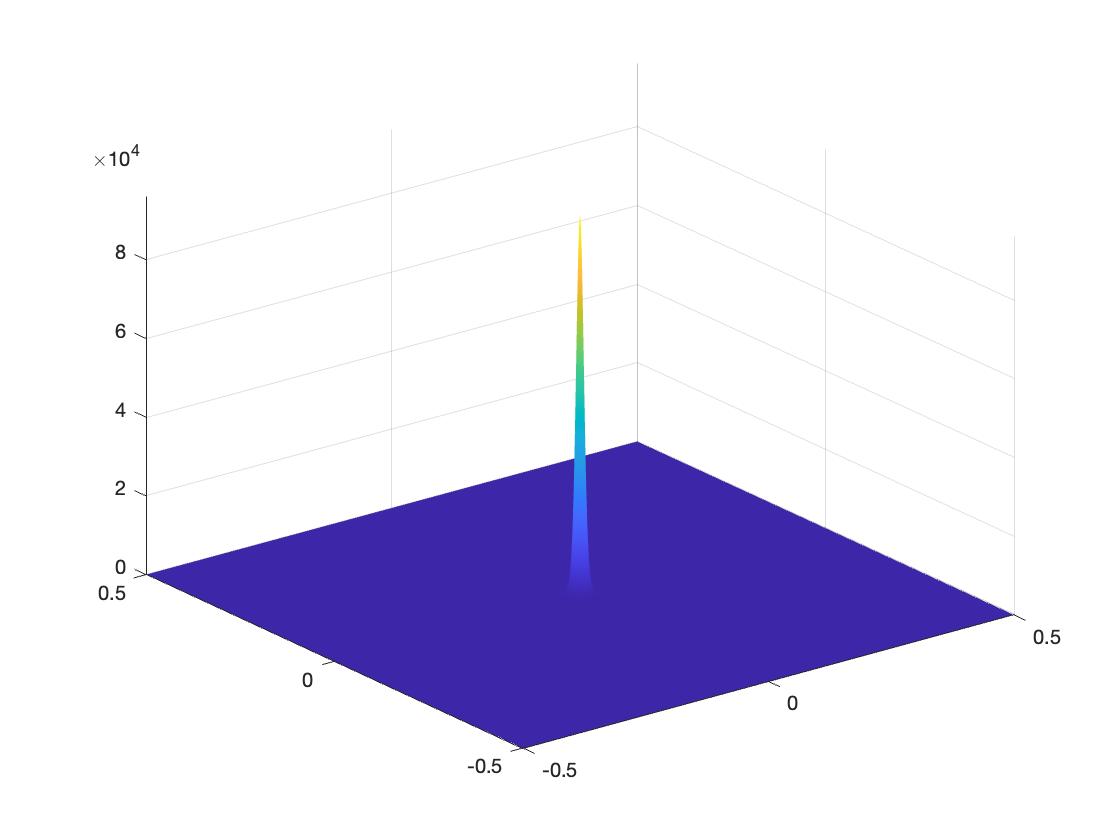}
    \caption*{$t = 1\times 10^{-4}$}
\end{minipage}
\hfill
\begin{minipage}{0.48\textwidth}
    \centering
    \includegraphics[width=\linewidth]{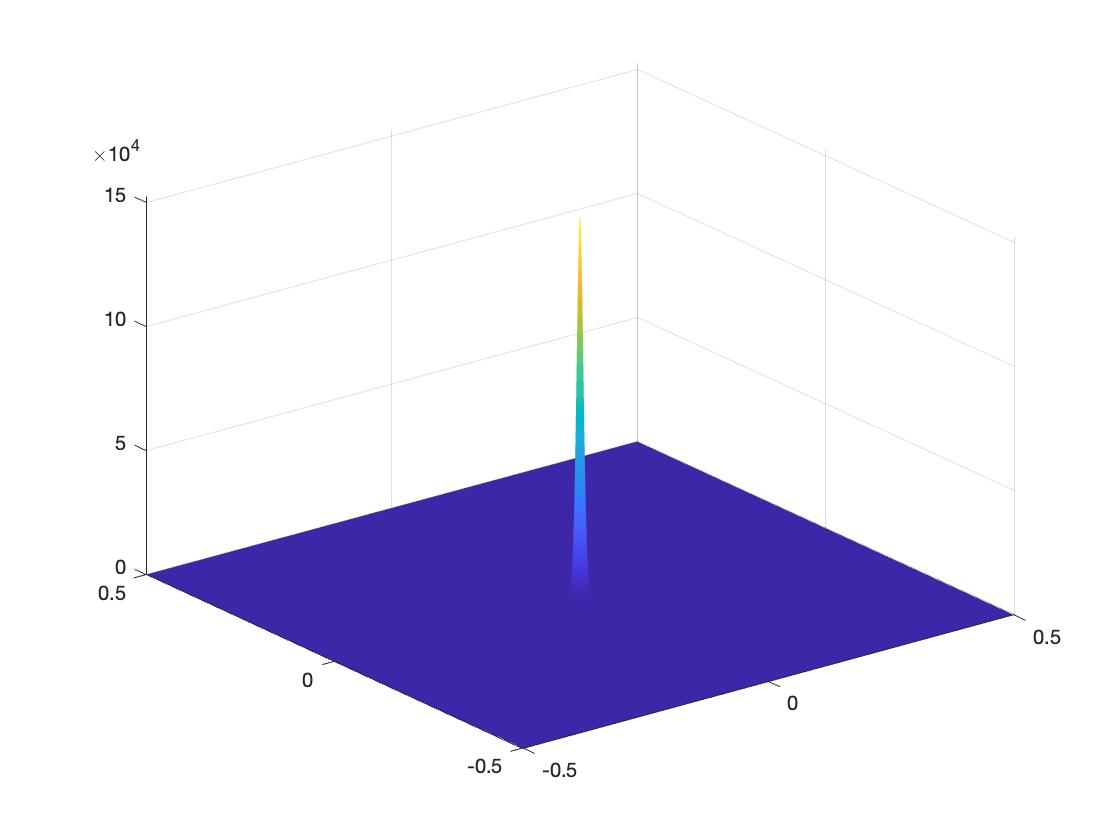}
    \caption*{$t = 2\times 10^{-4}$}
\end{minipage}

\caption{Evolution of \(\rho_1\) from \(t=3\times10^{-6}\) to \(t=2\times10^{-4}\).}
\label{fig:exampleD_blowup}
\end{figure}

We simulate this test case up to final time \(t = 2\times10^{-4}\) with time step \(\Delta t = 3\times10^{-6}\). The computed solution, shown in Figure~\ref{fig:exampleD_blowup}, exhibits a clear blow-up at the origin. The maximum of $\rho_1$ grows monotonically, reaching values of \(1.45\times10^5\) at the final time. No spurious negative values are observed, demonstrating the positivity-preserving property of the proposed method. The computed peak values are consistent with those reported in high-resolution simulations \cite{li2017local, QiuLiuYan2021, SulmanNguyen2019}.

To assess the sensitivity of the proposed method to the clustering exponent $\alpha$, we perform a parameter sweep over $\alpha \in \{3, 5, 6, 8, 10\}$, holding all other settings fixed: time step $\Delta t = 3\times10^{-6}$, sampling radius $r = 1$, and network configuration. For each value of $\alpha$, we record the resulting maximum density $\max \rho_1$ obtained with both the standard feedforward tanh architecture and the proposed MLP architecture (Fourier features + gated residual connections). Results are summarized in Table~\ref{tab:alpha_sensitivity}.

\begin{table}[htbp]
\centering
\caption{Peak density $\max \rho_1$ versus clustering exponent $\alpha$ in Example~D.}
\begin{tabular}{|c|c|c|}
\hline
$\alpha$ & Feedforward tanh ($\max \rho_1$) & MLP ($\max \rho_1$) \\
\hline
3  & $6.21\times 10^4$ & $4.30\times 10^4$ \\
5  & $8.10\times 10^4$ & $1.39\times 10^5$ \\
6  & $1.02\times 10^5$ & $1.45\times 10^5$ \\
8  & $1.10\times 10^5$ & $1.39\times 10^5$ \\
10 & $1.03\times 10^5$ & $1.40\times 10^5$ \\
\hline
\end{tabular}
\label{tab:alpha_sensitivity}
\end{table}

At $\alpha = 3$, the clustering is too weak to concentrate collocation points near the blow-up location, and both architectures underresolve the peak; in this regime the MLP offers no benefit over the feedforward baseline and in fact underperforms it by approximately 31\%, indicating that the added representational capacity of the Fourier-feature encoding cannot compensate for a poorly localized training set. For $\alpha \geq 5$, the proposed MLP architecture consistently and substantially outperforms the feedforward-tanh baseline, with gains ranging from approximately 26\% at $\alpha = 8$ to 72\% at $\alpha = 5$. The peak density is attained at $\alpha = 6$, where the MLP reaches $\max \rho_1 = 1.45\times10^5$, an improvement of approximately 42\% over the feedforward result at the same exponent. 

\subsubsection{Example E: Two-species parabolic-elliptic blow-up}

We consider the parabolic-elliptic limit (\(\varepsilon=0\)) of the two-species system on the domain \(\Omega=[-1,1]^2\). The parameters are
\[
\chi_1=1,\quad \chi_2=20,\quad \mu_1=\mu_2=1,\quad
D=1,\quad \alpha_1=\alpha_2=1,\quad \beta=1,
\]
and the chemoattractant satisfies the elliptic equation
\(-\Delta c = \rho_1+\rho_2 - c\) in \(\Omega\) with homogeneous Neumann boundary conditions.

The initial data are radially symmetric and identical for both species, given by the sharply peaked Gaussian profile
\[
\rho_1(x,y,0)=\rho_2(x,y,0)=50\,e^{-100(x^2+y^2)},
\]
which yields total masses \(M_1=M_2\approx\pi/2\). According to the criteria in Section~\ref{blow}, the blow-up condition is not satisfied, while the global existence criterion fails because \(M_2\) exceeds its critical threshold. Consequently, the problem lies in the undetermined regime, for which the theory guarantees neither global existence nor finite-time blow-up.

To address this open question, we rely on numerical evidence. High-resolution finite-volume simulations on uniform grids \cite{KurganovLukacovaMedvidova2014} show that \(\rho_2\) blows up, while \(\rho_1\) grows much more slowly. This behavior was later confirmed using an adaptive moving mesh method \cite{Chertock2019}, which captures the sharp spike in \(\rho_2\) with far fewer cells and demonstrates that both species blow up simultaneously, albeit with very different rates.

We replicate this challenging test case using our discrete-time PINN framework. The simulation is run up to final time \(t=4\times10^{-3}\) with time step \(\Delta t = 8\times10^{-5}\).

\begin{figure}[t]
\centering
\begin{minipage}{0.48\textwidth}
    \centering
    \includegraphics[width=\textwidth]{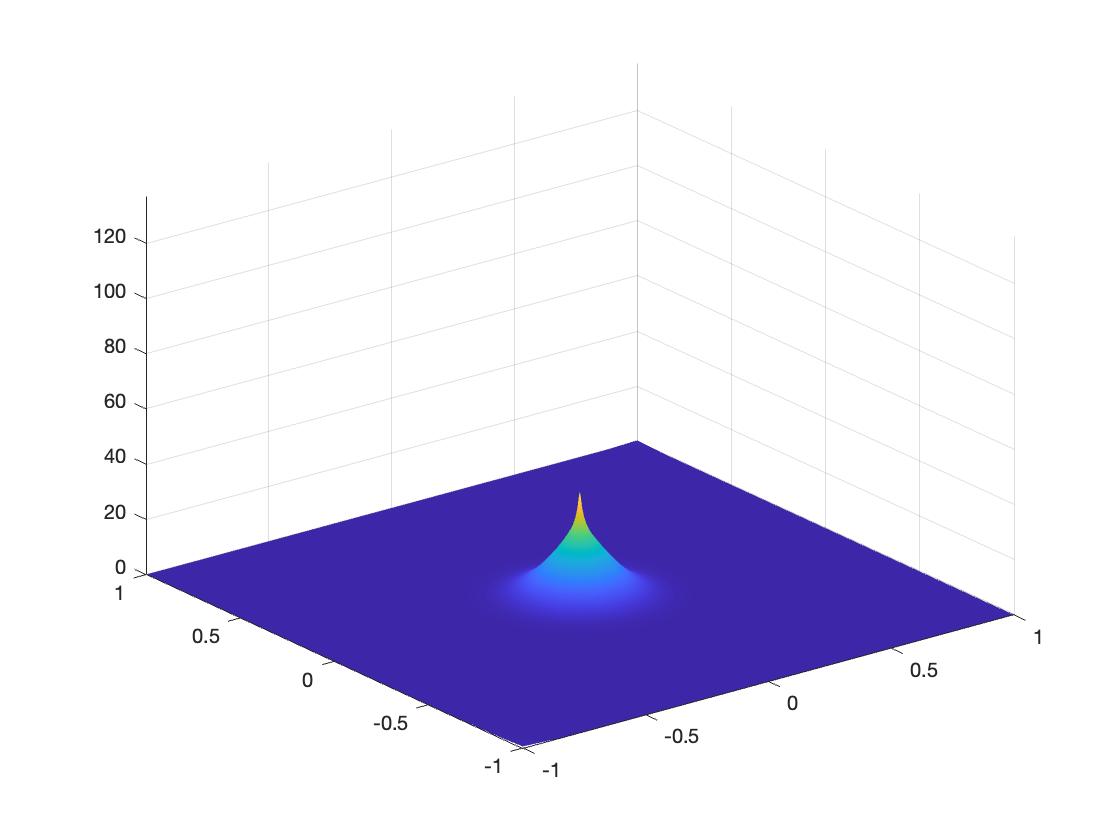}
    \caption*{$\rho_1$}
\end{minipage}
\hfill
\begin{minipage}{0.48\textwidth}
    \centering
    \includegraphics[width=\textwidth]{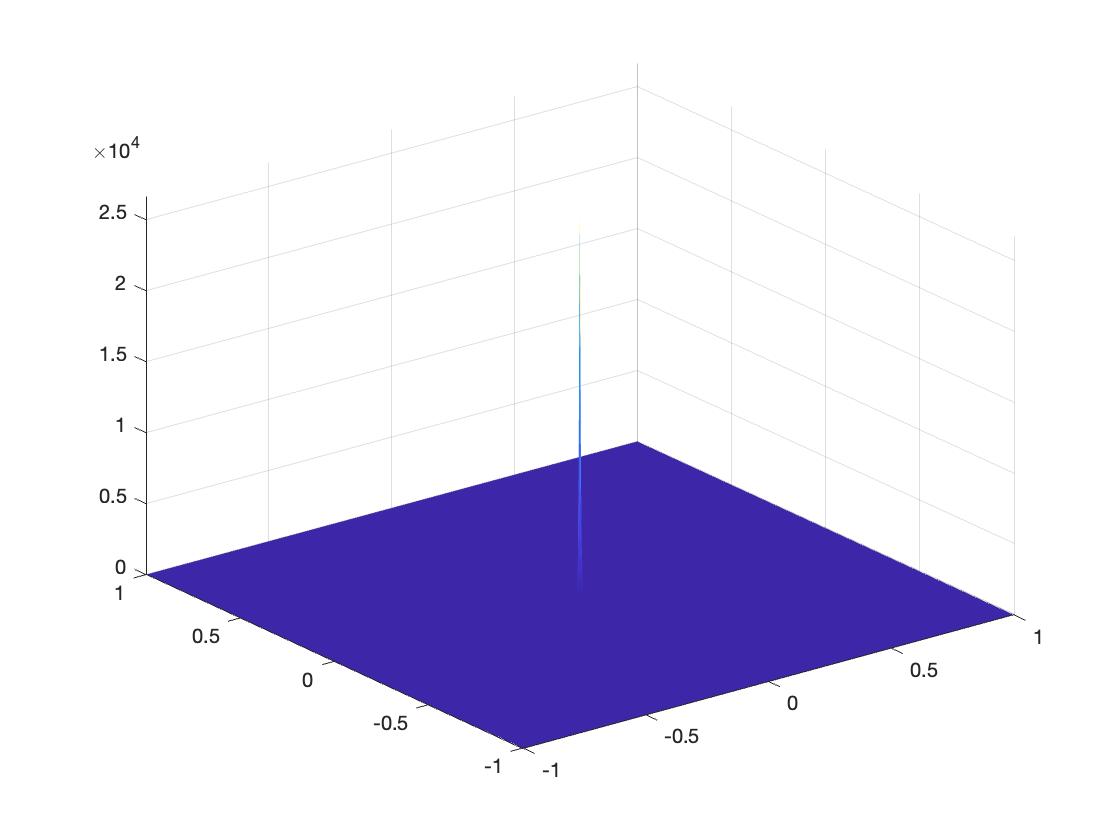}
    \caption*{$\rho_2$}
\end{minipage}
\caption{Time evolution and final profiles of $\rho_1$ and $\rho_2$ at \(t=4\times10^{-3}\).}
\label{fig:exampleE}
\end{figure}

Figure~\ref{fig:exampleE} shows that our PINN results agree qualitatively with the reference simulations: \(\rho_2\) develops a sharp, highly localized peak with \(\max \rho_2 \approx 2.87\times 10^4\), while \(\rho_1\) grows more moderately with \(\max \rho_1 \approx 3.72\times 10^2\). This behavior is consistent with the numerical results reported in~\cite{Chertock2019}, where \(\rho_2\) concentrates and blows up more rapidly than \(\rho_1\). The discrete-time PINN captures both components accurately without excessive dissipation, confirming its ability to resolve disparate blow-up scales.

\subsubsection{Example F: Two-species parabolic--parabolic blow-up}

We consider the fully parabolic Keller--Segel system ($\varepsilon=1$) with parameters
\[
\chi_1=5,\quad
\chi_2=60,\quad
\mu_1=\mu_2=1,\quad
D=10,\quad
\alpha_1=\alpha_2=1,\quad
\beta=1.
\]

First, we take the initial conditions on the domain \(\Omega=[-1,1]^2\):
\[
\rho_1(x,y,0)=\rho_2(x,y,0)
=500\,e^{-100(x^2+y^2)},
\qquad
c(x,y,0)=1.
\]
While a complete theoretical characterization of finite-time blow-up remains unavailable for the fully parabolic system, prior numerical studies have demonstrated that sufficiently large chemotactic sensitivities can trigger highly localized aggregation. Given the substantially stronger chemotactic sensitivity of the second species (\(\chi_2=60\)), it is expected to aggregate more rapidly and attain a higher peak concentration.

\begin{figure}[t]
\centering
\begin{minipage}{0.48\textwidth}
    \centering
    \includegraphics[width=\textwidth]{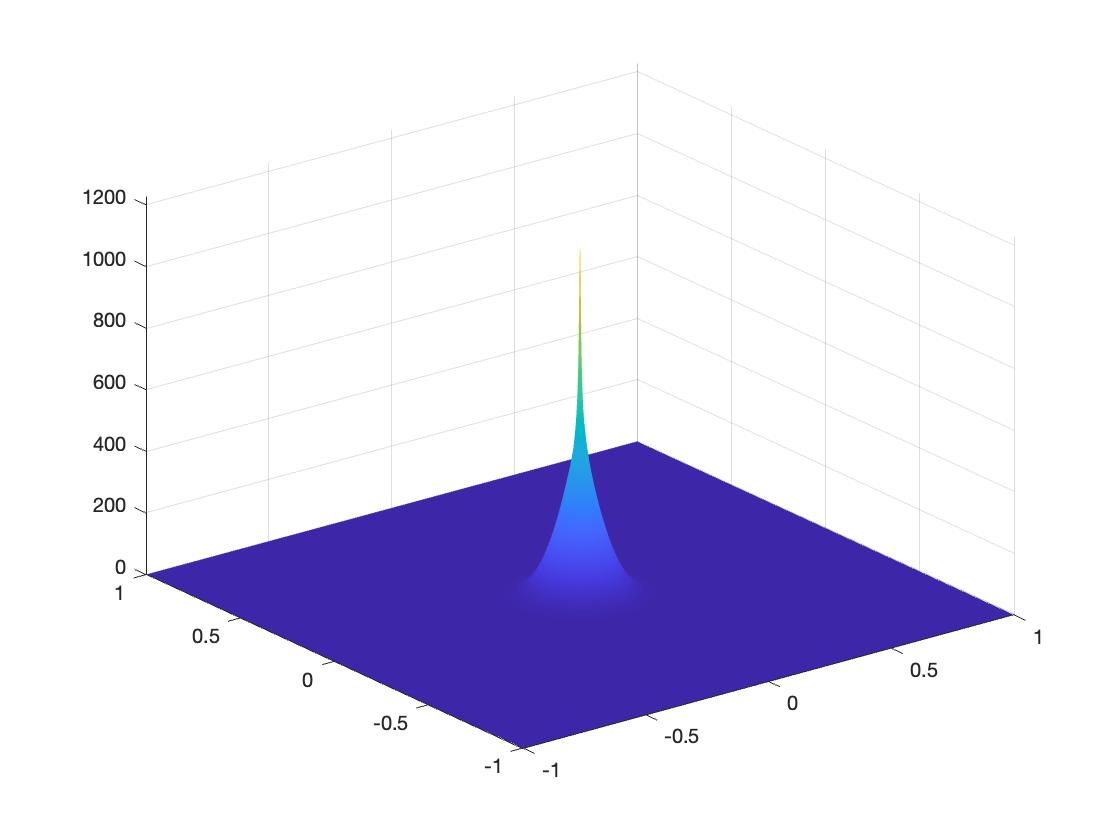}
    \caption*{$\rho_1$}
\end{minipage}
\hfill
\begin{minipage}{0.48\textwidth}
    \centering
    \includegraphics[width=\textwidth]{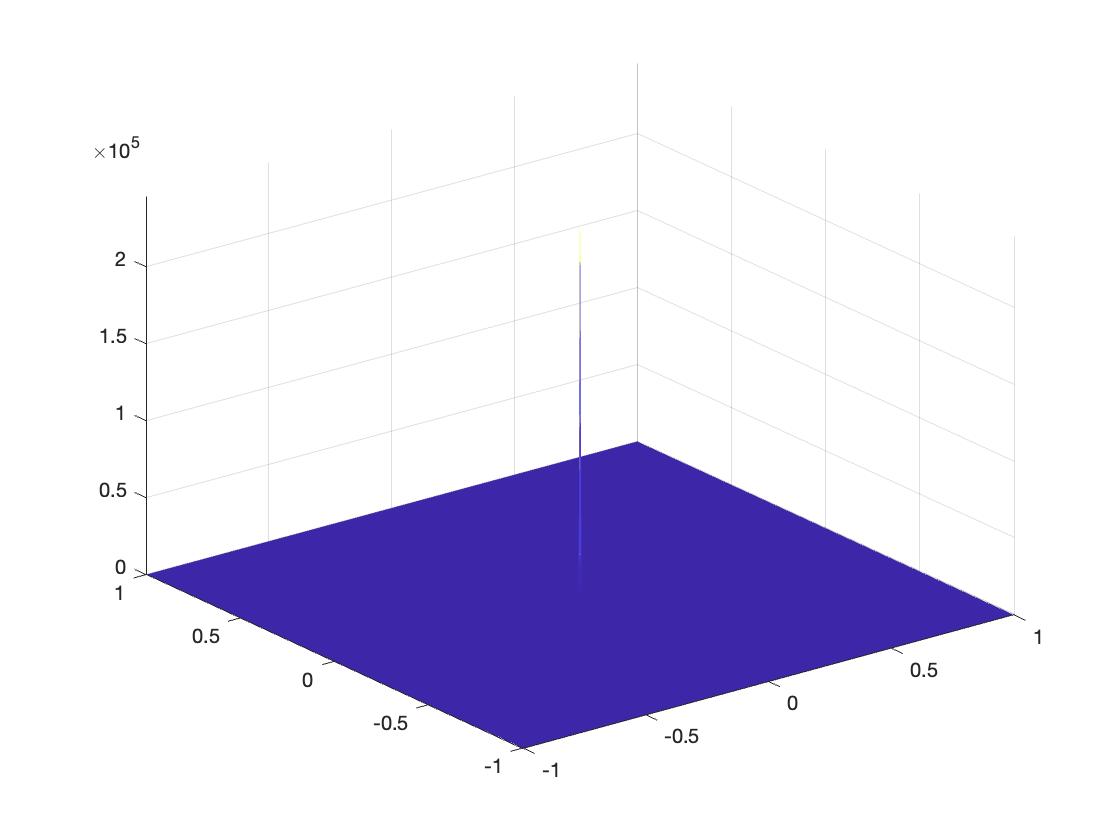}
    \caption*{$\rho_2$}
\end{minipage}
\caption{Blow-up profiles on \([-1,1]^2\) at \(t=8.4\times10^{-4}\).}
\label{fig:exampleF_lowmass}
\end{figure}

The simulation is run up to final time \(t=8.4\times10^{-4}\) with time step \(\Delta t=4\times10^{-5}\). Consequently, Figure~\ref{fig:exampleF_lowmass} shows that both populations have formed sharply localized peaks at the origin, with \(\rho_2\) reaching a substantially higher concentration (\(\max \rho_2 \approx 6.63\times 10^5\)) than \(\rho_1\) (\(\max \rho_1 \approx 1.37\times 10^3\)), consistent with its larger chemotactic sensitivity.

\begin{figure}[t]
\centering
\begin{minipage}{0.48\textwidth}
    \centering
    \includegraphics[width=\textwidth]{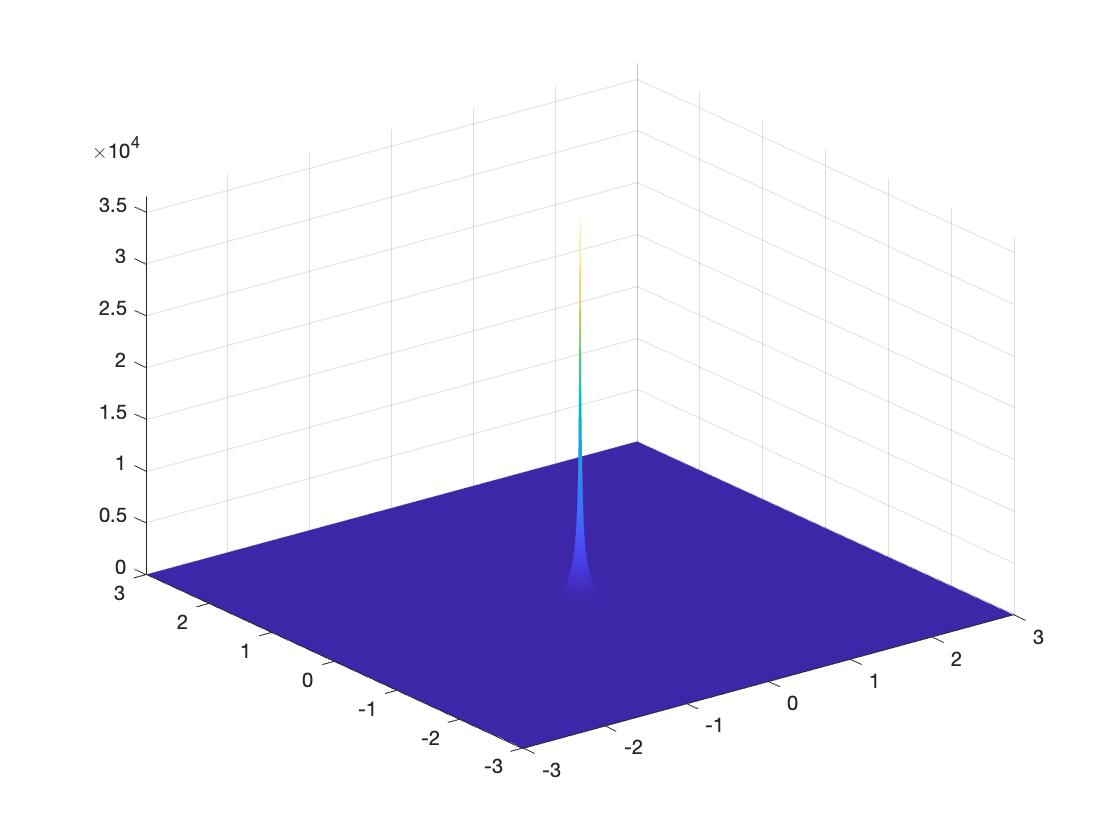}
    \caption*{$\rho_1$}
\end{minipage}
\hfill
\begin{minipage}{0.48\textwidth}
    \centering
    \includegraphics[width=\textwidth]{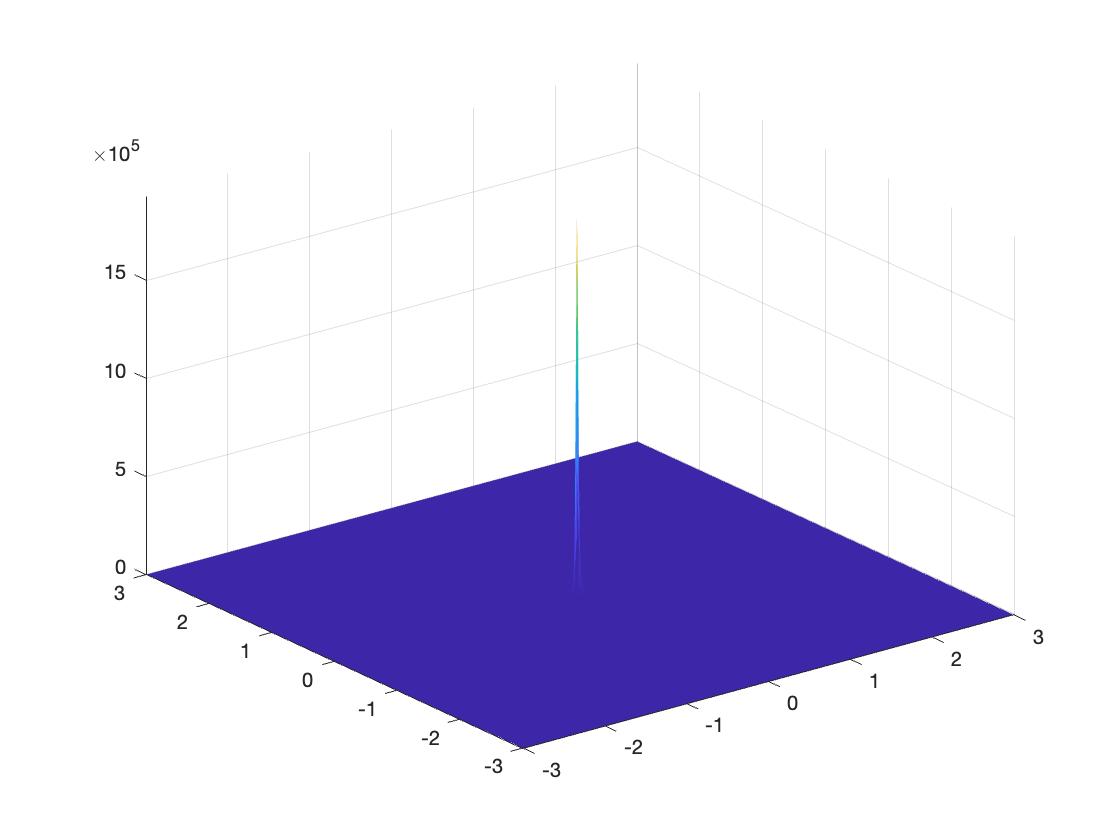}
    \caption*{$\rho_2$}
\end{minipage}
\caption{Blow-up profiles on \([-3,3]^2\) at \(t=1.65\times10^{-4}\).}
\label{fig:exampleF_highmass}
\end{figure}

To examine the effect of initial mass on the blow-up dynamics, we retain all model parameters but increase the amplitude of the initial densities by one order of magnitude on the larger domain \(\Omega=[-3,3]^2\):
\[
\rho_1(x,y,0)=\rho_2(x,y,0)
=5000\,e^{-100(x^2+y^2)},
\qquad
c(x,y,0)=1.
\]
For this more challenging case, we run up to final time \(t=1.65\times10^{-4}\) with time step \(\Delta t=5\times10^{-6}\). As shown in Figure~\ref{fig:exampleF_highmass}, the higher initial mass dramatically accelerates aggregation: \(\rho_2\) reaches \(\max \rho_2 \approx 1.98\times 10^6\), while \(\rho_1\) attains \(\max \rho_1 \approx 3.65\times 10^4\). Both species develop sharply localized profiles within a much shorter time, with \(\rho_2\) again dominating due to its stronger sensitivity. The blow-up onset advances from approximately \(t=8.4\times10^{-4}\) in the first experiment to \(t=1.65\times10^{-4}\) in the second. These outcomes are consistent with previously reported numerical experiments~\cite{KurganovLukacovaMedvidova2014, Chertock2019, li2017local}, demonstrating the robustness of the proposed method in capturing increasingly severe blow-up dynamics.

\section{Conclusion}
We have introduced a hybrid physics-informed neural network framework for multi-species chemotaxis that adapts its formulation to the qualitative regime of the underlying solution. In the subcritical regime, a standard continuous-time PINN with alternating training achieves relative errors on the order of \(10^{-3}\)--\(10^{-4}\) across one-dimensional, two-dimensional single-species, and fully coupled multi-species benchmarks, comparing favorably against established finite-element, finite-volume, and moving-mesh references as it can be seen in Tables~\ref{tab:1d_errors}--\ref{tab:multi_species_errors}. In the supercritical regime, the discrete-time formulation, combining backward Euler time stepping, a logarithmic output transformation, Fourier feature encoding with gated residual connections, and residual-adaptive collocation clustering, resolves the highly localized aggregation peaks and finite-time blow-up characteristic of single- and two-species Keller--Segel dynamics (Examples~D--F), reproducing peak magnitudes and growth rates consistent with high-resolution finite-volume and adaptive-mesh studies from the literature \cite{Chertock2019, KurganovLukacovaMedvidova2014, li2017local, QiuLiuYan2021, SulmanNguyen2019}. Both formulations preserve positivity intrinsically and operate mesh-free, resolving singular structures without post-processing limiters.

The central design principle is economy: rather than forcing one architecture to handle both smooth and singular dynamics, the added cost of discrete time-stepping, log-transformation, and adaptive clustering is incurred only where the solution demands it. The smooth-clustering strategy of Section~3.2 is essential to this: residual-based adaptive sampling alone introduces geometric artifacts that under-resolve the peak (Figure~\ref{fig:smooth_vs_irregular}), whereas the proposed symmetric clustering captures blow-up amplitudes accurately.

Two limitations define the scope of the present work and motivate future extensions. First, while network architecture hyperparameters were selected automatically using Optuna~\cite{Akiba2019Optuna}, extending this approach to the discrete-time scheme's collocation and discretization parameters, the clustering exponent \(\alpha\), sampling radius \(r\), and time step \(\Delta t\), proved computationally prohibitive, particularly in the supercritical regime where each training run is already expensive. These parameters were therefore tuned manually through the empirical sensitivity study in Table~\ref{tab:alpha_sensitivity}. Developing automated strategies for their selection, together with adaptive time-step control, would improve both robustness and efficiency. Second, because ground-truth solutions are unavailable near blow-up, supercritical validation relies on comparisons of peak magnitudes and growth trends against high-resolution reference simulations rather than rigorous error norms. The adaptive collocation strategy generalizes naturally to multiple aggregation centers, but a systematic demonstration of multi-spike cases awaits benchmarks with clearly identifiable, PINN-learnable blow-up behavior.

Future work will pursue extending automated hyperparameter search to the collocation parameters \(\alpha\), \(r\), and \(\Delta t\); adaptive time-step control; and extension of the hybrid framework to other PDEs exhibiting mixed smooth/singular behavior, such as reaction--diffusion systems with shock formation or finite-time extinction.



\end{document}